\documentclass[11pt,a4paper,Color]{article}
\usepackage{amssymb}
\usepackage{amsmath}
\usepackage{graphicx, color}
\usepackage{indentfirst}
\definecolor{bomi}{rgb}{0.4, 0.2, 0.8}

\input xypic.sty

\def\qed{\hfill {\large ${\sqcup\!\!\!\!\sqcap}$}}

\newenvironment{demo}{{\bf Proof }}
{\qed \\}

\newcommand{\re}{\mathbb R}

\newcommand{\eme}{{\rm M}^{n+1}_\lambda}
\newcommand{\flecha}{\longrightarrow}
\newcommand{\<}{\left<}
\renewcommand{\(}{\left(}
\newcommand{\lb}{\label}
\newcommand{\nn}{\nonumber}
\newcommand{\fracc}{\displaystyle\frac}

\renewcommand{\>}{\right>}
\renewcommand{\)}{\right)}

\newcommand{\us}{\underset}

\newcommand{\bde}{\begin{defi}}
\newcommand{\ede}{\end{defi}}

\numberwithin{equation}{section}
\newcommand\bal{\begin{align}}
\newcommand\eal{\end{align}}
\def\be{\begin{equation}}
\def\ee{\end{equation}}
\newcommand{\ba}{\begin{array}}
\newcommand{\ea}{\end{array}}

\def\oo{{\overline \omega}}
\def\oO{{\overline \Omega}}
\def\og{{\overline g}}
\def\oR{{\overline R}}
\def\oRic{{\overline Ric}}
\def\oSec{{\overline Sec}}
\def\oM{{\overline M}}

\def\oN{{\overline \nabla}}
\def\oD{{\overline \Delta}}
\def\ot{{\overline \theta}}

\def\oF{{\overline F}}
\def\t1{{\widehat {\mathfrak 1} }}
\def\tg{{\widehat g }}
\def\tR{{\widehat R }}
\def\tRic{{\widehat Ric}}
\def\tSec{{\widehat Sec}}
\def\tN{{\widehat \nabla}}
\def\tD{{\widehat \Delta}}
\def\tB{{\widehat {\rm B}}}

\def\ev{{\epsilon\nu}}
\def\te{{\widehat e }}
\def\tf{{\widehat f }}
\def\ttr{{\widehat r }}
\def\tt{{\widehat \theta }}
\def\tom{{\widehat \omega }}
\def\tOm{{\widehat \Omega }}

\def\tX{{\widehat X}}

\def\g{{\frak{g}}}
\def\a{\alpha}
\def\p{\varphi}
\def\r{\rho}
\def\fpp{{\fracc{\p'}{\p}}}
\def\fppp{{\fracc{\p''}{\p}}}
\def\s{{\rm s_\lambda}}
\def\c{{\rm c_\lambda}}
\def\sm{{\rm s_\mu}}
\def\cm{{\rm c_\mu}}
\def\sk{{\rm s_k}}
\def\ck{{\rm c_k}}

\def\ral{\sqrt{|\lambda|}}
\def\tr{{\rm tr}}

\def\dist{{\rm dist}}

\def\nablaa{\overline{\nabla}}
\def\Deltaa{\overline{\Delta}}
\def\vle{{\rm vol}}
\def\parcial#1#2{\frac{\partial #1}{\partial#2}}

\def\Daparcial#1#2{\frac{\nablaa #1}{\partial#2}}

\def\flecha{\longrightarrow}

\def\ds{\displaystyle}

\newtheorem{defi}{Definition}
\newtheorem{teor}{Theorem}

\newtheorem{prop}[teor]{Proposition}
\newtheorem{notacion}[defi]{Notation}
\newtheorem{lema}[teor]{Lemma}

\newtheorem{nota}{Remark}

\numberwithin{lemap}{teor}
\numberwithin{corop}{teor}

\numberwithin{ejer}{subsection}
\numberwithin{ejemplo}{subsection}

\begin{document}

\title{Mean curvature flow of graphs in warped products}



\author{ Alexander A. Borisenko and Vicente Miquel}


\date{ }

\maketitle

\begin{abstract}
Let $M$ be a complete Riemannian manifold which either is compact or has a pole, and let $\p$ be a positive smooth function on $M$. In the warped product $M\times_\p\re$, we study the flow by the mean curvature of a locally Lipschitz continuous graph on $M$ and prove that the flow exists for all time and that the evolving hypersurface is $C^\infty$ for $t>0$ and is a graph for all $t$. Moreover, under certain  conditions, the flow has a well defined limit.
\end{abstract}

\section{Introduction}

Let $M$ be a $n$-dimensional manifold, $(\oM,\og)$ a $n+1$ dimensional Riemannian manifold. A map $F:M \times [0,T[ \flecha \overline M$ such that every $F_t:= F(\cdot, t ) : M \flecha \oM$ is an immersion  is called the mean curvature flow (MCF for short) of $F_0$ if it is a solution of  the equation
\begin{align}\label{vpmf}
\parcial{F}{t} &=  \vec{H} 
\end{align}
where $\vec{H}(\cdot, t)$ is the mean curvature vector of the immersion $F_t$.

We shall use
the following {\it convention signs for the mean curvature $H$, the Weingarten map $A$ and the second fundamental form ($h$ for the scalar valued version and $\a$ for its vector valued version)}. For a chosen unit normal vector $N$, they are:

$A X = - \nablaa_XN$, $\a(X,Y) = \<\nablaa_X Y, N\> N = \<A X, Y\> N$, $h(X,Y)= \<\a(X Y), N\>$ and $H= \tr A = \sum_{i=1}^n h(E_i,E_i)$, for a local orthonormal frame $E_1, ..., E_n$  of the submanfold, and $\vec{H} = \sum_{i=1}^n \a(E_i,E_i) = H\ N$  . 

Along the rest of the paper, by 
$M_t$ we shall denote both the immersion $F_t:M\flecha \overline M$ and the image $F_t(M)$, as well as the Riemannian manifold $(M,g_t)$
 with the metric $g_t$ induced by the immersion. Analogous notation will be used when we have a single immersion $F:M\flecha \overline M$. Notice that, since in this paper we shall deal with graphs, all the immersions that will appear evolving by mean curvature flow will be, in fact, embeddings.
 
 In two fundamental papers \cite{EcHu0} and \cite{EcHu}, Ecker and Huisken studied the evolution of noncompact hypersurfaces in the Euclidean Space. In \cite{EcHu0} they studied the evolution of a   graph in $\re^{n+1}$ and showed that 
 
 (A) If $F_0$ a \lq\lq locally Lipschitz" continuous graph and has linear growth rate for its height, then \eqref{vpmf} with initial condition $F_0$ has a smooth solution for all time.
 
 (B) If, moreover,  $F_0$ is \lq\lq stright'' at the infinity, $F_t$ asymptotically approaches a selfsimilar solution of \eqref{vpmf}. They give an example showing that the condition cannot be weakened.
 
 In  \cite{EcHu}, Ecker and Huisken obtained some interior estimates and applied them to  prove  that the hypothesis of linear growth in (A) is not necessary, that is:
 
 (A') If $F_0$ a \lq\lq locally Lipschitz" continuous graph, then \eqref{vpmf} with initial condition $F_0$ has a smooth solution for all time.
 
 In \cite{Unt} and \cite{Un}, Unterberger extended result (A')  to the Hyperbolic Space $\mathcal H^{n+1}$ and gave a result of type (B). In this space, the first problem to face with  is the choosing of the right concept of \lq\lq graph''. A natural one is to say that a hypersurface $M$ of $\mathcal H^{n+1}$ is a graph over a totally geodesic hypersurface $\mathcal H^n$ if all the geodesics orthogonal to $\mathcal H^n$ cut $M$ once and transversally (we shall call it a {\it geodesic graph}). But, in his thesis \cite{Unt}, Unterberger found an example of hypersurface  which is a geodesic graph  but loses this property when it evolves under \eqref{vpmf}. Then he considered another concept of graph. Let $p\in \mathcal H^n$ be a fixed point, $\Gamma$ a geodesic through $p$ orthogonal to $\mathcal H^n$. We shall call {\it equidistant} curves all the curves which are at constant distance from $\Gamma$. Then we say that $M$ is an {\it equidistant graph}  over $\mathcal H^n$ if it cuts once and transversally all the equidistant curves.  Unterberger proved the exact analog of (A') for equidistant graphs. As a result of type (B), he proved that if $F_0: M \flecha \mathcal H^{n+1}$ is a \lq\lq locally Lipschitz"  equidistant graph and is at bounded distance from $\mathcal H^n$, then it converges asymptotically to $\mathcal H^n$. 
 
 After $\mathcal H^{n+1}$, other ambient spaces natural for trying to extend the results of Ecker and Huisken are products $M\times \re$ or, more generally, warped products $\re \times_\p M$ or $M\times_\p\re$. As usual, by a warped product $\mathcal M \times_\p \mathcal N$ of two Riemannian manifolds $(\mathcal M, g)$ and $(\mathcal N, h)$ we understand the riemannian manifold $(\mathcal M\times \mathcal N,  g + \p^2 h)$, being $\p:\mathcal M\flecha \re$ a positive smooth map. For these spaces it is natural to say that {\it a hypersurface of $M\times_\p\re$ (or  $\re \times_\p M$) is a graph } if it is a graph of a function $u:M\flecha \re$ and $u(M)$ cuts transversally the curves $s\mapsto (x,s)$. The interest in the last years for studying minimal and constant mean curvature surfaces in these ambient spaces  produces also a natural interest on the study of MCF on them.

 When $\p(s)= \cosh s$ and $M=\mathcal H^n$,  $\re\times_\p M$ is $\mathcal H^{n+1}$ and the concept of being a graph coincides  with that of geodesic graph in  $\mathcal H^{n+1}$. Then the counterexample of Unterberger gives us few hope of getting general results. Computations of the evolution of the gradient of $u$ (see the appendix) give some analytic reasons of the failing of the preservation of geodesic graphs under MCF. Also the last paragraph in Remark \ref{notagrafo} and the pictures in the appendix can help to get some gometric insight on this fact.
 
 When  $M=\mathcal H^n$, $x_0\in M$ and $\p(x)= \cosh (dist(x_0,x))$,  $M\times_\p \re$ is $\mathcal H^{n+1}$ and the concept of being a graph coincides with  that of equidistant graph in $\mathcal H^{n+1}$. Then,  general   $M\times_\p \re$ seem to be good general ambient spaces where to extend results of type (A') and (B). In this paper we show that this is in fact the case, proving an extension of (A') for general $M\times_\p \re$ (theorems \ref{existco},  \ref{stenonco} and \ref{steL}), where the only conditions to be satisfied by $M$ and $\p$ are: the quotients  $\fracc{|\tN^m\p|}{\p}$  are bounded, the curvature of $M$ and all its covariant derivatives are bounded and, when $M$ is non-compact, $M$ {\it has a pole} (that is, a point with empty cut-locus). May be the last condition can be weakened with some stronger analytic tools. On the other hand, like in \cite{EcHu}, the only condition for $M_0$ is to be a \lq\lq locally Lipschitz" graph. We also get a result of type (B) (theorems 	\ref{convcomp} and \ref{convnonco}) imposing on $M_0$ the condition that its distance to $M$ is bounded on it and on $M$ a pinching condition on its sectional curvature related to the pinching of the hessian of $\p$.
 
 We use mainly the methods in \cite{EcHu} and \cite{Un}, where we have to introduce necessary technical tricks to choose the right functions to get estimates having into account the complications introduced by the terms containing $\p$ and the curvature of $M$. We also  needed to use the comparison theory of Riemannian Geometry to bound some functions of the distance to a point or to a hypersurface. Moreover we cannot use barriers as it is done in \cite{Un} because we are not in a model space where we know the evolution of some hypersurfaces and, at some points, we need to substitute barriers arguments for others (see the end of the proofs of Theorem \ref{existnonco} and Lemma \ref{ltozero}). 
 
 Somewhat surprising for us has been the fact that the qualitative results of type (A') (not the estimates) do not depend on the curvature of $M$. 
 
 In contrast, the curvature of $M$ plays an essential role in the theorem of convergence.
  
 The paper is organized as follows: in section 2 we state the notation and recall some lemmas that will be used later. In section 3 we collect the properties of the ambient spaces and their hypersurfaces that we shall need.  Section \ref{evoluT} is a short comment about short time existence when $M$ is compact and also a useful description of the evolution of some geometric quantities under an equivalent flow in the direction of the last coordinate $u$ in $M\times_\p\re$. In section \ref{graph} we give the gradient etimate, which, when $M$ is compact, is enough to conclude  the preservation of the property of being a graph under MCF .  In section \ref{LTE} we obtain the higher order estimates which give rise to the long time existence (and conclude -module the existence theorems below- with the main theorems of the paper when the initial condition is smooth). In section \ref{STE} we complete the discussion of the existence when $M$ is non-compact, and, in section \ref{STEL} we discuss the existence theorems for Lipschitz initial conditions. It is more usual to give the theorems for smooth and Lipschitz initial conditions  simultaneously but, for the sake of non expert readers,  we preferred to do it separately.  In this way appears more clearly why in the Lipschitz case we cannot use the initial conditions to bound  the second derivatives and beyond, which forces us to introduce bounds depending on $t$ although these are becoming worse when $t$ goes to $0$.  Finally,  in section \ref{CONV} we give a case where  the flow has a limit and determine the limit.

{\bf Acknowledgments:} This work was done while the first author was Visiting
Professor at the University of Valencia in 2008,
supported by a \lq\lq ayuda del Ministerio de Educaci\'on y Ciencia SAB2006-0073.'' He wants to thank that university and
its Department of Geometry and Topology by the facilities they gave him. 

Second author  was partially supported  by DGI(Spain) and FEDER Project MTM2007-65852.

\section{ Preliminaries}
In this paper we shall consider a Riemannian manifold $(M,\tg)$ and an immersion $F: M\flecha \oM$ into another Riemannian manifold $(\oM, \og)$, and we shall denote by $g$ the metric  induced on $M$ by the immersion $F$. 

 We shall use  the notation $|X|$ or $\<X,Y\>$ for indicating the $\og$-norm, $\tg$-norm or $g$-norm of $X$ or the $\og$-product, $\tg$-product or $g$-product of $X$ and $Y$  if $X$ and $Y$ are, respectively,  tangent to $\oM$,  $M$ or $F(M)$. 
 
 For any vector $X\in T_{F(x)}M$, we shall denote by $X^\top$ the component of $X$ tangent to $F(M)$.

We shall use $\oN$, $\tN$ and  $\nabla$ to denote the covariant derivative and the gradient in $(\oM,\og)$,    $(M,\tg)$ and $(F(M), g)$ respectively. By $\oD$, $\tD$ and $\Delta$ we shall denote the corresponding laplacians.
The convention sign for the laplacians wiil be:
$\Deltaa f =  \tr \nablaa^2 f.$\\
For any  $\lambda\in\re$, we shall use the  notation: 
$$\s(t) \text{ is the solution of the equation $s'' + \lambda s=0$ satisfying $s(0)=0$ },  \qquad\c(t)=\s'(t) $$
These functions  satisfy the following computational rules:
\begin{equation}\label{trirel}
c_\lambda^2 + \lambda\ s_\lambda^2=1, \quad c_{4\lambda} =
c_\lambda^2 - \lambda\ s_\lambda^2, \quad s_{4\lambda}=
s_\lambda c_\lambda.
\end{equation}
When $\lambda<0$ one has $\s(t) = \sinh(\ral\  t)\big/\ral$ \ and \  $\c(t) = \cosh(\ral\  t)$.

Given a point $x_0$ in $M$, we
shall denote by $\ttr$ the  $\tg$-distance to
$x_0$ in $M$.  We shall denote by $\partial_{\ttr}^\bot$ at $x\in M$ the vector space orthogonal to $\partial_\ttr:= \tN \ttr$ in $T_xM$. Given a function $f: \re \flecha \re$, $f(\ttr)$ will mean $f\circ \ttr$.
The comparison theorem for the hessian of the distance function $\ttr$ says
\begin{lema}[cfr \cite{GW} or \cite{Pe}] If $\tSec \ge \mu$, then, at points of $M$ between $x_0$ and its cut«locus, one has 
\be \label{comphess}
\tN^2 \ttr \left|_{\partial_\ttr^\bot}\right.  \le \fracc{\cm}{\sm}(\ttr) \  \tg, \qquad \tD \ttr \le (n-1) \   \fracc{\cm}{\sm}(\ttr).
\ee
\end{lema}
 This will be essential for the gradient estimates we shall obtain in Theorem \ref{graphbo2nonco}. The next comparison theorem will be used for proving the convergence of the solution of \eqref{vpmf} in certain cases.
 \begin{lema}[cf. \cite{Gr}]
Let $\oM$ be a complete  riemannian manifold, with sectional curvature  $ \oSec \le k <0$, let $M$ be a complete totally geodesic hypersurface of $\oM$. Let   $\ell$ be the $\og$-distance to $M$, $\partial_\ell := \oN\ell$. Let us denote by $S_\ell$ the Weingarten map (associated to $\partial_\ell$)nof the hypersurface $\tau_\ell M$ at distance $\ell$ of $M$ (that is, $S_\ell X = -\oN_X\partial_\ell$ for every $X\in T\tau_\ell M$). At every point between $M$ and the cut locus of $M$ in $\oM$ one has 
\begin{align}
\<S_\ell X, X\> \le k \fracc{\sk}{\ck} |X|^2 \text{ for every }X\in T\tau_\ell M. \label{comptg}    
\end{align}
\end{lema}
Also in the proof of convergence we shall use the following maximum principle for noncompact manifolds due to Ecker and Huisken:

\begin{lema}[\cite{EcHu} Th. 4.3] \lb{echu43} Let $M$ be a manifold with a family of Riemannian metrics $g_t$ satisfying, for some $x_0$, 
$$\vle_{g_t}(B^{g_t}_r(x_0)) \le e^{k(1+r^2)},$$
where $B^{g_t}_r(x_0)$ is the $g_t$-ball of radius $r$ centered at  $x_0$,
Let $f: M\times [0,T] \flecha \re$ be a function which is smooth on $M\times ]0,T]$ and continuous on $M\times [0,T]$, satisfying:
\begin{itemize}
\item [(i)] $\ds \parcial{f}{t} \le \Delta^{g_t} f + \<\vec{a},\nabla^{g_t} f\> + b\ f$,  with $|\vec{a}|$ and $|b|$ bounded on $M\times [0,T]$
\item[(ii)] $f(\cdot,0) \le 0$
\item[(iii)] $\ds \int_0^T\(\int_M e^{\a r^2} |\nabla^{g_t} f|^2 d\mu_t\)dt < \infty$ for some $\a>0$
\item[(iv)] $ \sup_{M\times[0,T]} \left| \fracc{dg}{dt}\right| \le \beta$ for some $\beta <\infty$.
\end{itemize}
Then $f\le 0$ on $M\times[0,T]$.

\end{lema} 

\section{The geometric setting}

\subsection {The ambient space}

Given $(M,\tg)$ a $n$-dimensional riemannian manifold
 and $\varphi: M \flecha \re^+$ a $C^\infty function$, our ambient space will be
$\oM = (M\times \re, \og)$ with $\og = \tg +  \varphi(x)^2 \ du^2  $ (usually denoted as $M\times_{\p}\re$).

If $\tX$ (resp.$\tt$) is a vector field (resp. a differential form) on $M$, we shall denote by the same letter the vector fields (resp. forms) on $\oM$ induced by the family of embeddings  $j_s:M\flecha M\times\{s\} \subset M\times\re$ is defined by $j_s(x)= (x,s)$.

In a $\tg$-orthonormasl local frame $\te_1, ..., \te_n$  of $(M,\tg)$, if
$\tt^1, ..., \tt^n$ is the dual frame,  $\tg = \tt^{1^2} + ... + \tt^{n^2}$.
In the extension of this frame to $\oM$ completed with $du$, $\og$ is written as 
$\og = \varphi^2 du^2 +  (\tt^1)^2 + ... +  (\tt^n)^2$.
Then, if we denote $\ot^0 = \p  du,\  \ot^1=\ \ \tt^1, \ ...\ , \ot^n = \ \tt^n$, one has
$$
\og = (\ot^0)^2  + (\ot^1)^2 + ... +\ ...\ , +(\ot^n)^2
$$
and its dual frame
${\overline e}_0, \  {\overline e}_1, ..., {\overline e}_n$  in $\oM$
 is a $\og$-orthonormal frame related with the $\te_i$ 
  by ${\overline e}_i = \te_i$,
  and ${\overline e}_0=  \ds\frac{1}{\varphi} \parcial{}{u} =:  \ds\frac{1}{\varphi} \partial_u$

\medskip
\noindent From now on $i,j = 1, ..., n$ and $a,b= 0,1,...,n$. In this subsection we shall use he Einstein convention of summing repeated indices when one is a subindex and the other a super-index.

\subsubsection{The Levi-Civita connection of $\oM$}

In the frames given above, we compute the Cartan connection forms $\oo_a^b$ defined by
 \be d\ot^b = - \sum_{a=0}^n\oo^b_a \wedge  \ot^a \label{Cartan1}\ee
Differentiating the formulae for $\ot^a$ and comparing with \eqref{Cartan1}, we have
\begin{align}
d \ot^0 &= d\p\wedge du  = \te_i(\p) \tt^i \wedge \fracc1{\p} \ot^0= -  \fracc1{\p} \te_i(\p) \ot^0  \wedge \ot^i = -\oo_j^0 \wedge \ot^j \nn \\
d \ot^i &=  d\tt^i =  -  \tom_j^i \wedge \tt^j  = - \oo_j^i \wedge \ot^j- \oo^i_0 \wedge \ot^0 \nn
\end{align}
and the solution of these two equations is, denoting $\p_i:=\te_i(\p)$, 
\be
\oo_0^i = - \frac{\p_i}{\p} \  \ot^0 = -\p_i \ du, \quad \oo_j^i= \tom_j^i. \lb{oob}
\ee
The relation between $\oo$ and $\nablaa$ is 
$$\nablaa_XY = X(Y^a) {\overline e}_a + Y^a \oo_a^b(X) {\overline e}_b,$$
then we have
\be
\ba{ll}
 \oN_{{\overline e}_0} {\overline e}_0 = \oo_0^i({\overline e}_0) {\overline e}_i = - \fracc{1}{\p} \tN \p    & \oN_{{\overline e}_0} {\overline e}_i = \oo_i^j({\overline e}_0) {\overline e}_j+\oo_i^0({\overline e}_0) {\overline e}_0 = \fracc{\p_i}{\p} {\overline e}_0 =\fracc{\p_i}{\p^2} \partial_u  \\
\oN_{{\overline e}_i} {\overline e}_0 = \oo_0^j({\overline e}_i) {\overline e}_j =0   & \oN_{{\overline e}_i} {\overline e}_j= \oo_j^k({\overline e}_i) {\overline e}_k + \oo^0_j({\overline e}_i) {\overline e}_0 =  \tN_{\te_i}\te_j .  \label{codeb}
\ea
\ee

\subsubsection{The curvature of $\oM$} \lb{curvoM}

In the frame ${\overline e}_a$, the curvature tensor $\oR$ and the curvature $2$-forms $\oO$ are related by
\be
\oR_{abcd} = \oO_c^d({\overline e}_a,{\overline e}_b) \nn
\ee
and the Cartan equations say
\be\lb{Cartan2}
d\oo_b^a = - \oo^a_c \wedge \oo^c_b - \oO_b^a.
\ee

Differentiating in \eqref{oob} and applying equations \eqref{Cartan2}, 
\begin{align} \label{fcurvb}
\oO_j^i = \tOm_j^i - \oo^i_0 \wedge \oo^0_j = \tOm_j^i  ,  \qquad
 \oO_0^i =  \fracc{ \p_{ij} }{\p} \ot^j \wedge \ot^0 - \fracc{ \p_j}{\p} \tom^i_j \wedge \ot^0. 
\end{align}
where $\tOm^i_j$ are the curvature $2$-forms  of $(M,\tg)$.

From the expressions for the curvature $2$-forms we obtain for the components of the curvature tensor:

\begin{align}
\oR_{0ijk} &= \oO_j^k({\overline e}_0, {\overline e}_i) = \tOm_j^k ({\overline e}_0,{\overline e}_i) = 0 \nn \\
\oR_{0i0k} &=  \oO_0^k({\overline e}_0, {\overline e}_i) = \( \fracc{ \p_{kj} }{\p} \ot^j \wedge \ot^0 - \fracc{ \p_j}{\p} \tom^k_j \wedge \ot^0\) ({\overline e}_0,{\overline e}_i) = -\fracc{\p_{ki} }{\p}  + \frac{\p_j}{\p} \tt^j\(\tN_{\te_i} \te_k\) = -\fracc{\tN^2 \p}{\p}(\te_i,\te_k) \nn \\
\oR_{ijk\ell} &=  \oO_k^\ell({\overline e}_i, {\overline e}_j) =  \tOm_k^\ell ({\overline e}_i,{\overline e}_j) =  \tR_{ijkl} . \lb{curoM}
\end{align}

\begin{nota}\label{bounR}
From  the formulae for $\oo$ and $\oR$ given above it follows that the norms (in the metric induced by $\og$) of $ \oR$ and $\oN^i \oR$ ($i=0,1,2, ... $) are bounded if the norms (in the metric induced by $\tg$) of $\tR$, $\tN\tR$,   $\fracc{\tN^i\p}{\p}$, ($i=1,2, ... $)   are bounded. 

The {\bf general setting} {\it in this paper is that the geometry of $\oM$ is bounded, that is, $|\oR|$, $|\oN^i \oR|$ ($i=0,1,2, ... $) are bounded.}
\end{nota}

\begin{nota}\label{fh0}
From the formulae for the covariant derivatives it follows that

The hypersurfaces $u=constant$ are totally geodesic in $\oM$. 

Fixed a $x_0\in M$, if $\tN\p(x_0)=0$,  the curve $u\mapsto(x_0,u)$ is a geodesic and  the curves $u\mapsto (x,u)$ have constant geodesic curvature $\fracc{|\tN\p|}{\p}$, and all the points in $u\mapsto (x,u)$ are at constant distance from the curve $u\mapsto (x_0,u)$ (are equidistant curves).

\end{nota}

Now we use the expression  of the curvature components  to obtain the  sectional curvatures:
\begin{align}
&\overline S_{i0} = \oR_{i0i0} = -\fracc{\tN^2 \p}{\p}(\te_i,\te_i) , 
\qquad \overline S_{ij} = \oR_{ijij}=   \widehat S_{ij} ,  \label{scoMb} 
\end{align}

\begin{nota} When $M$ is the simply connected space of constant sectional curvature $\lambda$, $ M_\lambda^n$ and  $\p(x) = \c(\ttr(x))$ ($\ttr: M \flecha \re$ defined by $\ttr(x) = \dist_{\tg}(x_0,x)$ for some $x_0$ fixed), then $\oM$ must be $\eme$), the simply connected space of constant sectional curvature $\lambda$, as can be checked from the above formulae for curvatures.
\end{nota}


\subsection{The submanifold}\label{STEb}

\subsubsection{The vectors and quantities $N$, $\nabla u$,  $\sigma$ and $v$}

We consider an embedding $F: M\flecha \oM$ given by the graph $F(x) =(x,u(x))$ of a function $u: M \flecha \re$.

The frame $\te_i$ of $M$ induces a frame (on the submanifold $F(M)$)  $e_i = F_*\te_i =  \te_i + u_i \partial_u = {\overline e}_i + u_i \p {\overline e}_0$, where $u_i =\te_i(u)$. In this frame, the matrix of the metric $g$ of the submanifold and its inverse, and the dual frame $\theta^i$ are given by:
\be
\ba{l}
g_{ij} = \p^2 \ u_i\ u_j +  \delta_{ij},\\ 
g^{ij} = \delta^{ij} - \fracc{\p^2 u_i u_j}{1 + \p^2 |\tN u|^2} \lb{g^ij} 
\ea
\ee

A unit normal vector to $F(M)$ can be found using $\xi = a\tN u + b \partial_u$, imposing the condition $\<\xi,e_i\>=0$ and dividing by the $\og$-norm. We choose
 \begin{align*}
 N&= \fracc{- \p^2 \tN u +  \partial_u}{\p \sqrt{\p^2 |\tN u|^2 + 1}} = \fracc{- \p \tN u +  {\overline e}_0}{ \sqrt{\p^2 |\tN u|^2 + 1}}.\\
 \text{ Then } &\<N,\partial_u\> =  \fracc{ \p}{ \sqrt{\p^2 |\tN u|^2 + 1}} \text{ and } \<N,{\overline e}_0\> =  \fracc{ 1}{ \sqrt{\p^2 |\tN u|^2 + 1}} 
 \end{align*}
 
 The gradient of $u$ in $F(M)$ can be computed using 
  \be \oN u= \fracc{{\overline e}_0}{\p}, \quad \nabla u = \oN u^\top = \fracc1{\p}\({\overline e}_0 - \<{\overline e}_0, N\> N\) =: \fracc1\p {\overline e}_0^\top\lb{gradu}  
\ee
  or $\nabla u = g^{ij} u_i e_j$. In both cases we obtain:
\be \label{fgradub}
\nabla u = \fracc{\tN u + |\tN u|^2 \partial_u}{1+\p^2  |\tN u|^2}
\ee
From \eqref{gradu} or \eqref{fgradub} one gets
\be
|\nabla u|^2 = \fracc1{\p^2}\(1 -\< N, {\overline e}_0 \>^2\) =  \fracc{|\tN u|^2}{1+ \p^2 |\tN u|^2 }
\ee
and
\be
|\tN u|^2 = \fracc{ |\nabla u|^2}{1 - \p^2 |\nabla u|^2} = \fracc{1 -\<N,{\overline e}_0\>^2}{\p^2 \<N,{\overline e}_0\>^2} 
\ee

If we define $\sigma = \<N,{\overline e}_0\>$ and $v = \fracc1{\sigma}$, the above formulae read
\be\lb{nuv}
  |\nabla u|^2 =\fracc1{ \p^2} \(1 - \fracc1{v^2}\), \qquad |\tN u|^2 = \fracc{ v^2-1}{\p^2}.
\ee

From the above formulae we also get
$$N = \fracc{- \p\ \tN u + e_0}{v}, \qquad \tN u = \fracc{- v N + {\overline e}_0}{\p} \qquad \tN u = \fracc{ \p^2 v^2 \nabla u -  (v^2 -1) \partial_u.}{\p^2}$$

\subsubsection{Relations of $\nabla v$ and $\nabla^2 u$ with the second fundamental form}

For computations, on $F(M)$ we shall consider another local frame $E_i$ orthonormal and satisfying $(\nabla_{E_i}E_j)_p=0$ at the point $p$ where we are doing the computations. Using it we compute the hessian of $u$. First:
$$\oN_{E_i}{\overline e}_0 = \<E_i, {\overline e}_0\> \oN_{{\overline e}_0}{\overline e}_0 = - \<E_i, {\overline e}_0\>  \fracc{\tN\p}{\p}.$$
\begin{align}
(\nabla^2 u)(E_i,E_j) & = \< \nabla_{E_i} \nabla u, E_j\> =  \< \nabla_{E_i} \fracc1{\p}({\overline e}_0 -\<{\overline e}_0, N\> N), E_j\>\nn \\
&= - \fracc{1}{\p} \<E_i,{\overline e}_0 \>  \<\fracc{\nabla\p}{\p}, E_j\> - \fracc{1}{\p} \<E_j,{\overline e}_0 \>  \<\fracc{\nabla\p}{\p}, E_i\> + \<{\overline e}_0, N\> \fracc1{\p} h(E_i,E_j) 
\end{align}
its trace
\begin{align}
\Delta u &= - \fracc2{\p^2} \<{\overline e}_0^\top,\nabla \p\> + \fracc1{\p v}   H     .\label{Dub}
\end{align}
and its norm
\begin{align}
|\nabla^2 u|^2 &= \fracc{4}{\p^4} \<{\overline e}_0^\top,\nabla \p\>^2+ \fracc1{\p^2} \sigma^2 |A|^2 - 4 \  \fracc1{\p^3} \  \sigma  \  h({\overline e}_0^\top,\nabla \p) ,
\end{align}
where ${\overline e}_0^\top := {\overline e}_0 - \<{\overline e}_0, N\> N = \p \nabla u$.

For $\nabla v$ we have, for every $X$ tangent to $F(M)$,
\begin{align*}
\<\nabla v,X\> &= -\fracc1{\sigma^2} \<\nabla \sigma, X\> = -\fracc1{\sigma^2} \(\<\oN_X{\overline e}_0,N\> + \<{\overline e}_0,\oN_X N\>\) \\
&= -\fracc1{\sigma^2} \(-\<X,{\overline e}_0\> \<\fracc{\tN\p}{\p},N\> - \<{\overline e}_0, AX \>\) 
= \fracc1{\sigma^2} \( \<  \<\fracc{\tN\p}{\p},N\> {\overline e}_0 + A{\overline e}_0^\top, X \>\) ,
\end{align*}
then
\begin{align}
\nabla v &= v^2 \(\<\fracc{\tN\p}{\p},N\> {\overline e}_0^\top + A{\overline e}_0^\top\) \lb{nablav} \end{align}
and
\begin{align}
\<\nabla v , \nabla \p\>&=  \<\nabla v , \oN \p\>\ = v^2 \< \<\fracc{\tN\p}{\p},N\> {\overline e}_0^\top+ A{\overline e}_0^\top,  \oN \p\> \nn\\
&= v^2 h({\overline e}_0^\top, \nabla \p) + \fracc{v^2}{\p} \<\tN\p, N\> \<\nabla\p,{\overline e}_0^\top\>.
\label{nvnp}
\end{align}

\section{Evolution of $u$ up to tangential diffeomorphisms and short time existence when $M$  is compact.} \lb{evoluT}

As it is well known, the equation \eqref{vpmf} is not parabolic, but it is also known (see\cite{Eck}, or also \cite{CaMi1})) that a solution $F(p,t)$ of \eqref{vpmf} can be written (as far as it remains a graph over $M$) as the composition $F(p,t) = \overline F(\Phi(p,t),t)$, where $\Phi(\cdot,t)$ is a family of diffeomorphismd depending on $t$ and $\overline F$ satisfies
  \be
\< \parcial{\overline F}{t}, N\> = H . \label{mcf}
\ee
and can be written, for each $t$, as a graph $\overline F(x,t) = (x, u(x,t))$ over $M$.  Using this parametrization of $\overline F$ the equation \eqref{mcf} becomes
\be \label{mcfg}
\<\partial_u , N\> \parcial{u}{t} = H, \quad \text{ that is } \quad  \parcial{u}{t} = \fracc{v}{\p} H 
\ee
In order to apply the standard theory of P.D.E., we compute $H$ in terms of $u:M\flecha \re$ and $\tN$. We shall use the vectors $e_i$ and the expression for $N$ given in section \ref{STEb}. The second fundamental form can be computed using the expressions \eqref{codeb} for $\oN$, which gives:
\begin{align}
\< A e_i , e_j\> & = -\< \oN_{e_i} N, e_j\> = \fracc{1}{\sqrt{\p^2  |\tN u|^2 + 1}} \< - \oN_{\te_i + u_i \partial_u} (- \p \tN u + {\overline e}_0), \te_j + u_j \p {\overline e}_0\> \nn \\
& = \fracc1v \( \p_i \ u_j + \p \<\tN_{\te_i}{\tN u} , \te_j\> + \p^2 u_i u_j\ \<\tN u , \tN\p\> + u_i   \p_j \)\lb{Aeiej}
\end{align}
From \eqref{Aeiej} and \eqref{g^ij}, 
\begin{align}
H = g^{ij} \< A e_i , e_j\> & =\(\delta^{ij} - \frac{\p^2 u_i u_j}{v^2}\)  \fracc1v \( \p_i \ u_j + \p \<\oN_{\te_i}{\tN u} , \te_j\> + \p^2 u_i u_j\ \<\tN u , \tN\p\> + u_i   \p_j \) \nn \\
 &= \fracc1v \(2\<\tN \p, \tN u\>  + \p \tD u + \p^2 |\tN u|^2 \<\tN u , \tN\p\> \) \nn \\
 & \qquad  - \fracc{\p^2}{v^3} \( 2 |\tN u|^2 \<\tN u,\tN \p\> + \p \< \tN_{\tN u} \tN u, \tN u\> + \p^2 |\tN u|^4  \<\tN u,\tN \p\> \) \nn \\
 & = \fracc{\p}{v} \( \tD u - \fracc{\p^2}{v^2}  \< \tN_{\tN u} \tN u, \tN u\> \)  + \fracc1v \<\tN u,\tN\p\>  \fracc{v^2+1}{v^2}
\lb{HI}
\end{align}
By substitution of this expression in \eqref{mcfg},
\begin{align}\label{mcft}
\parcial{u}{t} 
& = \tD u - \fracc{v^2-1}{v^2}  \tN^2_{\t1 \t1} u  + \fracc1\p \<\tN u,\tN\p\>  \fracc{v^2+1}{v^2}
\end{align}
which is a parabolic equation whereas $v$ remains bounded.

Then, when $M$ is compact, the existence of solution of \eqref{mcft} for a small interval of time  and smooth initial condition follows directly from the theory of parabolic equations. When $M$ is non-compact, we postpone the discussion to section \ref{STE}.

For sections \ref{STE} and \ref{STEL} it will be useful to know the evolution under \eqref{mcf} of $N$, $v$ and $H$. We get them now.

Using the parametrization $\oF(x,t)=(x, u(x,t))$ and \eqref{codeb},
\begin{align*}
\Daparcial{N}{t} &= \sum_{ij} g^{ij}\<\oN_{\parcial{\overline F}{t}} N, e_i \>  e_j = - \sum_{ij} g^{ij} 
 \< N, \oN_{\parcial{\overline F}{t}}  e_i \>  e_j = - {\parcial{u}{t}}  \p \sum_{ij} g^{ij} 
 \< N, \oN_{{\overline e}_0 } e_i  \>  e_j \nn\\
 & =- {\parcial{u}{t}} \sum_{ij} g^{ij} \< N, \oN_{{\overline e}_0 } (\te_i + u_i \p {\overline e}_0) \>  e_j = \parcial{u}{t}  \( - \< N, {\overline e}_0\>\fracc{\nabla\p}{\p} + \<N,\tN\p\> \nabla u \) \end{align*}
 From the above evolution equation and the definition of $v$,
 \begin{align}
\parcial{v}{t} &= - v^2 \<\Daparcial{N}{t}, {\overline e}_0\> = v^2 \parcial{u}{t} \( \< N, {\overline e}_0\> \<\fracc{\nabla\p}{\p}, {\overline e}_0\> - \<N,\tN\p\>  \<\nabla u,  {\overline e}_0\> \) , \text{ then }\nn
\end{align}
 \begin{align}
\left|\parcial{v}{t}\right| &\le  v^2 \left|\parcial{u}{t}\right| \( \left|\fracc{\tN\p}{\p}\right| + |\tN\p| \) \lb{bo_var_v}
\end{align}
From \eqref{HI} we get the following evolution for $H$

\begin{align}
\parcial{H}{t} &= - \fracc1{v^2} \({\p} \(
 \tD u - \fracc{\p^2}{v^2}  \< \tN_{\tN u} \tN u, \tN u\>  \) 
 +  \<\tN u,\tN\p\>  \fracc{v^2+1}{v^2} \) \parcial{v}{t}\nn\\
 &+ \fracc{\p}{v} \(\tD(\frac{v}{\p} H) + 2 \fracc{\p^2}{v^3} \< \tN_{\tN u} \tN u, \tN u\> \parcial{v}{t} \right. \nn \\
 &\qquad \quad \left.- \fracc{\p^2}{v^2}  \( \< \tN_{\tN \(\frac{v}{\p} H\)} \tN u, \tN u\>  +  \< \tN_{\tN u} \tN (\frac{v}{\p} H), \tN u\>  +  \< \tN_{\tN u} \tN u, \tN \(\frac{v}{\p} H\)\> \) \) \nn \\
 &+\fracc1{v} \(\<\tN\(\frac{v}{\p} H\), \tN\p\> \fracc{v^2+1}{v^2} - \<\tN u, \tN\p\> \frac{2}{v^3} \parcial{v}{t} \) \label{var_HtT}
\end{align}

The following observation will also be  useful

\begin{prop}\lb{propax}
Let $f$ be a $C^2$ function defined on $M_t$, $\tf = f \circ \oF(\cdot,t)$,  then  $|\tN\tf|$ is bounded if $|\nabla f|$ and $|\tN u|$ are bounded, and $|\tN^2 \tf|$ is bounded if $|\nabla^2f|$, $|\nabla f|$, $|A|$, $|\tN u|$ and $|\tN^2 u|$ are bounded
\end{prop}
\begin{demo}
Let  $\tX$ be a vector field over $M$, since $X:=\oF(\cdot,t)_*\tX = \tX+ \<\tN u, \tX\> {\overline e}_0$,
\begin{align}
 \<\tN \tf, \tX\> &=  df(\oF(\cdot,t)_*\tX)  = \< \nabla f, \oF(\cdot,t)_*\tX  \> = \< \nabla f, \tX + \<\tN u, \tX\> {\overline e}_0\> \nn \\
 &= \< \nabla f, \tX \>+\<\nabla f, {\overline e}_0\>  \<\tN u, \tX\> ,\quad \text{ then} \nn\\
\tN \tf &=\nabla f - \< \nabla f, {\overline e}_0 \> {\overline e}_0 +\<\nabla f, {\overline e}_0\>  \tN u \lb{tNtfNf}
\end{align}
Using the above expressions for $X$ and $\tN\tf$ and the formulae \eqref{codeb} for $\oN$,
\begin{align} \tN_\tX \tN\tf  & =  \oN_{X - \<\tX,\tN u\> {\overline e}_0}  \tN \tf =  \oN_{X}  \tN \tf  - \<\tX,\tN u\> \oN_{{\overline e}_0}  \tN \tf \nn \\
=& \nabla_X\nabla f + \a(X,\nabla f)- \< \nabla_X\nabla f + \a(X,\nabla f), {\overline e}_0 \> {\overline e}_0  \nn \\
& \  + \<\nabla f,\<\tX,\frac{\tN\p}{\p}\> {\overline e}_0 - \<\tN u,\tX\>\frac{\tN\p}{\p}\> (-{\overline e}_0+\tN u)+ \< \nabla_X\nabla f + \a(X,\nabla f), {\overline e}_0 \> \tN u   \nn \\
&\  + \<\nabla f, {\overline e}_0\>  \(\tN_X\tN u + \< \tX, \tN u\> \<\tN u,\frac{\tN\p}{\p} \> {\overline e}_0\)  - \<\tX,\tN u\> \<\tN\tf,\frac{\tN\p}{\p}\> {\overline e}_0 .\lb{tN2fN2f}
\end{align}
It follows from \eqref{tNtfNf} that $|\tN\tf|$ is bounded if $|\nabla f|$ and $|\tN u|$ are bounded, and from \eqref{tN2fN2f} that $|\tN^2 \tf|$ is bounded if $|\nabla^2f|$, $|\nabla f|$, $|A|$, $|\tN u|$ and $|\tN^2 u|$ are bounded.   
\end{demo}

\section{Preserving the property of being a graph}\label{graph}
The condition  \lq\lq cuts transversaly the curves $s\mapsto (x,s)$" in the definition of a graph means that $\sigma:=\< N, {\overline e}_0\> >0$, which is equivalent to say $1< v = \fracc{1}{\sigma}< \infty$. Therefore, our first goal is to obtain an upper bound for $v$. To achieve this,  we need  the evolution equation for $v$ under \eqref{vpmf}.

\begin{lema} Under \eqref{vpmf}, $v$ evolves as
 \begin{align}
\parcial{v}{t} & = \Delta v - \frac2{v} |\nabla v|^2  + \fracc2{\p} \<\nabla v, \nabla\p\>  - v  |A|^2 \nn \\
& \qquad    - v \(1-\fracc1{v^2}\) \(\frac{\tD \p}{\p}  + \tRic_{\t1 \t1} +  \fracc{|\tN \p|^2}{\p^2}-\fracc{\tN^2\p}\p  \(\t1,\t1\)\) 
 \label{dvtnb}
\end {align}
where \lq\lq$\t1$" is the unit vector in the direction $\tN u$.
\end{lema}
\begin{demo} First we compute $\Delta \sigma$. To do so, we shall use an orthonormal frame of $M$ of the form $ E_1, ... , E_n$ as was introduced before. 
\be\label{Deltaub}
\Delta \sigma  = E_i E_i \<N,{\overline e}_0\> = E_i\(\<\oN_{E_i}N, {\overline e}_0\> + \<N,\oN_{E_i}{\overline e}_0\>\)  
\ee
\begin{align}
E_i \<\oN_{E_i}N, {\overline e}_0\> &= -E_i(h(E_i, {\overline e}_0^\top)) = -\nabla_{E_i}(h)(E_i,{\overline e}_0^\top) - h(E_i, \nabla_{E_i} {\overline e}_0^\top) \nn\\
& = -\nabla_{{\overline e}_0^\top}(h)(E_i,E_i)+\oR_{E_i {\overline e}_0^\top E_i N} - \<AE_i, \oN_{E_i} ({\overline e}_0-\<{\overline e}_0, N\> N)\> \nn\\
& = -{\overline e}_0^\top(H)+\oRic ({\overline e}_0^\top, N) + \<AE_i, \<E_i,{\overline e}_0\> \fracc{\tN\p}{\p}\> - \<AE_i, \<{\overline e}_0, N\> AE_i\> \nn \\
  &= -{\overline e}_0^\top(H)+\oRic ({\overline e}_0^\top, N) + \<A {\overline e}_0^\top, \fracc{\tN\p}{\p}\> - \<{\overline e}_0, N\>  |A|^2 \lb{Ds1sb}
\end{align}
\begin{align}
E_i  \<N,\oN_{E_i}{\overline e}_0\> &= E_i \< N, -\<E_i,{\overline e}_0\> \fracc{\tN \p}{\p})\> \nn \\
 & = -{\overline e}_0^\top \< N, \fracc{\tN\p}{\p}\>  - \<N , \fracc{\tN\p}{\p}\> \( \< \oN_{E_i}E_i ,{\overline e}_0\> + \<E_i,\oN_{E_i}{\overline e}_o\>\) \nn \\
 &=  - {\overline e}_0^\top \< N, \fracc{\tN\p}{\p}\>  - \<N , \fracc{\tN\p}{\p}\> \(  H \<N ,{\overline e}_0\> - \<E_i, \<E_i,{\overline e}_0\>\fracc{\tN\p}{\p}\>\) \nn \\
 &=  - {\overline e}_0^\top \< N, \fracc{\tN\p}{\p}\>  - \<N , \fracc{\tN\p}{\p}\>   \<N ,{\overline e}_0\> H + \<N , \fracc{\tN\p}{\p}\> \< {\overline e}_0^\top ,\fracc{\tN\p}{\p}\>  .\lb{Ds2sb}
\end{align}
Then
 \begin{align}
\Delta\sigma & =  -{\overline e}_0^\top(H)+\oRic ({\overline e}_0^\top, N) + \<A {\overline e}_0^\top, \fracc{\tN\p}{\p}\> - \<{\overline e}_0, N\>  |A|^2 \nn\\
& - {\overline e}_0^\top \< N, \fracc{\tN\p}{\p}\>  - \<N , \fracc{\tN\p}{\p}\>   \<N ,{\overline e}_0\> H + \<N , \fracc{\tN\p}{\p}\> \< {\overline e}_0^\top ,\fracc{\tN\p}{\p}\>  .\lb{DsTb}
\end{align}
Now, let us compute 
\begin{align}
\oRic ({\overline e}_0^\top, N) &= \oRic({\overline e}_0, N) -\<{\overline e}_0,N\>  \oRic(N,N)  \nn \\
&=  \oRic\({\overline e}_0, \fracc{- \p\ \tN u + {\overline e}_0}{v}\) - \fracc1v  \oRic\(\fracc{- \p\ \tN u + {\overline e}_0}{v},\fracc{- \p\ \tN u + {\overline e}_0}{v}\) \nn \\
& = \fracc1v \oRic ({\overline e}_0, {\overline e}_0) - \fracc{\p}{v} \oRic({\overline e}_0,\tN u) - \fracc{\p^2}{v^3} \oRic(\tN u, \tN u) + 2 \fracc{\p}{v^3} \oRic({\overline e}_0, \tN u) -\fracc{1}{v^3} \oRic({\overline e}_0, {\overline e}_0) \nn \\
&= \fracc1v \(1-\fracc1{v^2}\) \oRic({\overline e}_0,{\overline e}_0) -\fracc{\p}{v} \(1-\frac{2}{v^2}\) \oRic({\overline e}_0,\tN u)- \fracc{\p^2}{v^3} \oRic(\tN u, \tN u).
\end{align}
From the expressions \eqref{curoM} of the components of $\oR$, it follows that
\begin{align}
\oRic({\overline e}_0,{\overline e}_0)& = \oR({\overline e}_0, {\overline e}_i,{\overline e}_0,{\overline e}_i) = - \frac{\tD \p}{\p} \\
\oRic({\overline e}_0,\tN u) &= 0 \\
\oRic(\tN u,\tN u) &= \tRic(\tN u,\tN u) =  \fracc{v^2-1}{\p^2} \tRic_{\t1\t1}.
\end{align}
 then
 \begin{align}
 \Delta\sigma & =  -{\overline e}_0^\top(H) - \fracc1v \(1-\fracc1{v^2}\) \frac{\tD \p}{\p}  - \fracc1v \(1-\fracc1{v^2}\)  \tRic_{\t1\t1}.+ \<A {\overline e}_0^\top, \fracc{\tN\p}{\p}\> - \fracc1v  |A|^2 \nn\\
& - {\overline e}_0^\top \< N, \fracc{\tN\p}{\p}\>  - \<N , \fracc{\tN\p}{\p}\>   \<N ,{\overline e}_0\> H + \<N , \fracc{\tN\p}{\p}\> \< {\overline e}_0^\top ,\fracc{\tN\p}{\p}\>  .\lb{DsTb2}
 \end{align}

 \begin{align}
 - {\overline e}_0^\top \< N, \fracc{\tN\p}{\p}\>  &= - \fracc1\p  \oN^2\p({\overline e}_0^\top,N) - \fracc1\p \<\oN\p,\oN_{{\overline e}_0^\top} N\> + \fracc1{\p^2} \<{\overline e}_0^\top, \oN\p\> \<N,\oN\p \> \nn \\
&= - \fracc{\oN^2\p}\p  ({\overline e}_0^\top,N) + \fracc1\p h({\overline e}^\top, \nabla \p) + \fracc1{\p^2} \<{\overline e}_0^\top, \nabla\p\> \<N,\oN\p \>
 \end{align}
 which gives
  \begin{align}
 \Delta\sigma & =  -{\overline e}_0^\top(H) - \fracc1v \(1-\fracc1{v^2}\) \frac{\tD \p}{\p}  - \fracc1v \(1-\fracc1{v^2}\)  \tRic_{\t1\t1}. - \fracc1v  |A|^2 +  \fracc2\p h({\overline e}^\top, \nabla \p) \nn\\
&- \fracc{\oN^2\p}\p  ({\overline e}_0^\top,N) + \fracc2{\p^2} \<{\overline e}_0^\top, \nabla\p\> \<N,\oN\p \>  - \<N , \fracc{\tN\p}{\p}\>   \<N ,{\overline e}_0\> H  \text{, that is, using \eqref{nvnp}}\nn
\end{align}
\begin{align}  
 \Delta\sigma& =  -{\overline e}_0^\top(H) - \fracc1v \(1-\fracc1{v^2}\) \frac{\tD \p}{\p}  - \fracc1v \(1-\fracc1{v^2}\)  \tRic_{\t1\t1} - \fracc{\oN^2\p}\p  ({\overline e}_0^\top,N)  - \fracc1v  |A|^2  \nn\\
&+ \fracc2{\p\ v^2} \<\nabla v, \nabla\p\>   - \<N , \fracc{\tN\p}{\p}\>   \fracc Hv .\lb{DsTb3}
 \end{align}
 
On the other hand:
\be \nabla v = -\fracc1{\sigma^2} \nabla \sigma \ee
\be
\Delta v = E_iE_i v = E_i\(-\fracc1{\sigma^2} E_i\sigma\) = \fracc2{\sigma^3} |\nabla\sigma|^2 -\fracc1{\sigma^2} E_iE_i(\sigma) = \frac2{v} |\nabla v|^2 - v^2 \Delta \sigma.
\ee
\begin{align}
\parcial{v}{t} &  = -\fracc1{\sigma^2} \parcial{\sigma}{t} = -v^2 \(\<\Daparcial{N}{t}, {\overline e}_0\>  + \<N,\Daparcial{{\overline e}_0}{t} \>\) \nn \\
&= -v^2 \(\<-\nabla H, {\overline e}_0\> -  \< N , H \<N,{\overline e}_0\> \fracc{\tN\p}{\p}\> \) \nn \\
&= v^2 {\overline e}_0^\top (H) + v H \<N,\fracc{\tN\p}{\p}\> .
\end{align}
Joining the expressions for $\Delta v$, $\Delta \sigma$ and $\ds\parcial{v}{t}$, we obtain
 \begin{align}
\parcial{v}{t} & = \Delta v - \frac2{v} |\nabla v|^2  - v \(1-\fracc1{v^2}\) \(\frac{\tD \p}{\p}  + \tRic_{\t1\t1}\)  - v^2 \fracc{\oN^2\p}\p  ({\overline e}_0^\top,N)  - v  |A|^2  + \fracc2{\p} \<\nabla v, \nabla\p\> \nn\\
&  \nn 
\end{align}
Defining  $N^h:=N-  \<N,{\overline e}_0\> {\overline e}_0 = N - \sigma {\overline e}_0$,
\begin{align}
\oN^2\p ({\overline e}_0^\top,N) &= \oN^2\p \({\overline e}_0 -  \sigma (N^h + \sigma {\overline e}_0) , N^h + \sigma {\overline e}_0\) \nn \\
&= \sigma \fracc{|\tN \p|^2}{\p} - \sigma  \tN^2 \p(N^h,N^h) - \sigma^3 \fracc{|\tN \p|^2}{\p}  \nn \\
&= \fracc1v \(1- \fracc1{v^2}\) \fracc{|\tN \p|^2}{\p} - \fracc1v \tN^2 \p(N^h,N^h).
\end{align}
and sustitution in the evolution equation for $v$ gives
 \begin{align}
\parcial{v}{t} & = \Delta v - \frac2{v} |\nabla v|^2  - v \(1-\fracc1{v^2}\) \(\frac{\tD \p}{\p}  + \tRic_{\t1 \t1} +  \fracc{|\tN \p|^2}{\p^2}-\fracc{\tN^2\p}\p  \(\fracc{N^h}{|N^h|},\fracc{N^h}{|N^h|}\)\) \nn \\
& \qquad    - v  |A|^2  + \fracc2{\p} \<\nabla v, \nabla\p\>
 \end{align}
 Let us remark that $N^h = \fracc{- \p \tN u}{ \sqrt{\p^2 |\tN u|^2 + 1}}$, then $\t1 =\fracc{N^h}{|N^h|}$, and $\fracc{\tD \p}{\p}  -\fracc{\tN^2\p}\p  \(\fracc{N^h}{|N^h|},\fracc{N^h}{|N^h|}\)$ is $-\oRic({\overline e}_0,{\overline e}_0) + \oR_{0\t1 0\t1}$.
 \end{demo}
\begin{notacion}
A general hypothesis in this work is that $\oM$ has bounded geometry, in particular that $\fracc{|\tN\p|}{\p}$, $\fracc{|\tN^2\p|}{\p}$ and $|\tR|$  are bounded. Then there are constants $\eta$, $\mu_1$, $\mu_2$ and $\mu$ satisfying
\begin{align}
\eta = \sup_{\oM}\fracc{|\tN\p|}{\p}, \quad \mu_1 \tg \le \fracc{\tN^2\p}{\p} \le \mu_2 \tg, \quad \tSec \ge \mu \lb{def_mueta}
\end{align}
and we define the constant $\nu$ by the formula
\be
\mu = \fracc{- n \mu_1 + \mu_2}{n-1} - \nu \lb{def_nu}
\ee
\end{notacion}

\begin{teor}\label{graphbo2comp} 
Let $M$ be compact. Let $F(x,t):M \flecha \oM$ be a solution of  \eqref{vpmf}  defined on a  time interval $[0,T[$. If $\oM = M \times_{\p}\re$ and $M_0$  is a graph over $M$, then $M_t $ is a graph over $M$ for every $t\in [0,T[$. With more precision: $v(F(x,t)\le \max_{M_0} v \ e^{(n-1) \nu t}$.
\end{teor}
\begin{demo}  Using the notation \eqref{def_mueta} and \eqref{def_nu}, we have 
\be
 - \( \tRic_{\t1\t1}   - \fracc{\tN^2\p}{\p}(\t1,\t1) +\fracc{\tD \p}{\p} +\fracc{|\tN \p|^2}{\p^2}\) \le (n-1) \nu. \lb{boundRic}
 \ee
  From  \eqref{dvtnb} and forgetting about the negative summands, we reach the inequality:
\begin{align}
\parcial{v}{t} & < \Delta v - \frac2{v} |\nabla v|^2 +  \fracc2\p  \<\nabla v , \nabla \p\> + (n-1) \nu \ v . \label{evolvIb}
\end{align}
 By the maximum principle, $ v$ is bounded by the solution of $y'(t)= (n-1) \nu \ y(t)$ with the initial condition $y(0)=v_0 := \sup_{M_0} v$,  then $ \ds v \le v_0 e^{(n-1) \nu t}$, as claimed.

On the other hand, if $M_t$ is a graph over a proper subset of $M$ and $M$ is compact, it cannot be homeomorphic to $M$, then if $M_0$ is a graph over all $M$ it  cannot move to $M_t$ without producing singularities, which means $t\notin [0,T[$. From this and the estimate on $v$ it follows that $M_t$ remains to be a graph aver $M$ for all time in $[0,T[$.
\end{demo}

\vspace{-3mm}
\begin{teor}\label{graphbo2nonco} 
Let $M$ be complete non compact with a pole $x_0$. Let $M_t$ be a solution of  \eqref{vpmf}  defined on a maximal time interval $[0,T[$. If $\oM = M \times_{\p}\re$  and $M_0$ is a graph over $M$, then  $v\circ F(x,t)$ remains bounded at each point $F(x,t)$ for every $t\in[0,T[$. 

{\bf Remark} The hypothesis $M$ non compact also imposes $\mu \le 0$ (cf. \eqref{def_mueta}). On the other hand, if $\tSec\ge 0$, it is also true that $\tSec>\mu'$ for every $\mu'<0$, then we can suppose, without loose of generality, that $\mu<0$.
\end{teor}
\vspace{-3mm}
\begin{demo}  Let us recall that  $\ttr$ denotes  the $\tg$-distance from $x_0$.  As in the compact case, we have the inequality \eqref{boundRic}. An standard procedure to apply the maximum principle to $v$ in the non-compact case is to look at the evolution of the product $\phi v$ of $v$ times some cut-function $\phi$ suitably chosen. Inspired in \cite{Eck} and \cite{Un}, we define first a function 
$$\zeta(\ttr,t) = f(\ttr) \ e^{\beta n t},$$ 
where the function $f$ and the constant $\beta$ will be defined later, and 
\begin{align}
\phi: [0,\infty[ \flecha \re, \text{ by } \phi(\zeta) = \a\ (f(\r) - \zeta)^2, 
\end{align}
where, again, the constant $\a >0$ will be chosen later. Let us observe that:
\be
\phi' = - 2 \a\ (f(\r)-\zeta) <0 \text{ if } \zeta < f(\r), \quad \text{ and } \quad
\phi'' = 2 \a >0
\ee
In order to compute the evolution of $\phi v$ we need first to compute the 

\noindent {\bf evolution of $\ttr$}. Since $\oN\ttr =  \tN\ttr$, which we shall denote $\partial_{\ttr}$,
\begin{align}
&\(\parcial{}{t} - \Delta\) \ttr = H N(\ttr) - \sum_{i=1}^n E_iE_i\ttr = H \<N,\partial_\ttr\> - \sum_{i=1}^n \<\oN_{E_i} E_i,\partial_\ttr\>  - \sum_{i=1}^n\<E_i, \oN_{E_i} \partial_\ttr\>
\nn \\ 
&\qquad  = -  \sum_{i=1}^n \<E_i, \oN_{E_i} \partial_\ttr\>  = - \sum_{i=1}^n \<E_i, \oN_{E_i-\<E_i, \partial_\ttr\> \partial_\ttr -  \<E_i, {\overline e}_0\> {\overline e}_0} \partial_\ttr\> -\sum_{i=1}^n \<E_i,{\overline e}_0\> \<E_i, \oN_{{\overline e}_0}\partial_\ttr\>. \nn
\end{align}
Denoting $E_i^\bot := E_i-\<E_i, \partial_\ttr\> \partial_\ttr - \<E_i, {\overline e}_0\> {\overline e}_0 $,
\begin{align}
\(\parcial{}{t} - \Delta\) \ttr & =- \sum_{i=1}^n \tN^2 \ttr \(E_i^\bot,E_i^\bot\) - \sum_{i=1}^n \<E_i,{\overline e}_0\> \fracc{\partial_\ttr\p}{\p} \<{\overline e}_0,E_i\> \nn \\
&\us{\text{by }\eqref{comphess}}{\ge} -\fracc{\cm}{\sm} \ \sum_{i=1}^n \tg(E_i^\bot, E_i^\bot) -  \fracc{\partial_\ttr\p}{\p} \<{\overline e}_0^\top,{\overline e}_0^\top\> \nn \\
&= -\fracc{\cm}{\sm} \(n -|\partial_\ttr|^2 - |{\overline e}_0^\top|^2\) - \fracc{\partial_\ttr\p}{\p} |{\overline e}_0^\top|^2. \lb{328}
\end{align}
{\bf evolution of $f(\ttr) := f\circ\ttr$} for $f:\re\flecha \re$, with $f'>0$
\begin{align}
\(\parcial{}{t} - \Delta\)f(\ttr) &= f' \(\parcial{}{t} - \Delta\) \ttr - f'' \left|\partial_\ttr^\top\right|^2 \nn \\
&\us{\text{by }\eqref{comphess}}{\ge} f' \(-\fracc{\cm}{\sm} \(n -|\partial_\ttr^\top|^2 - |{\overline e}_0^\top|^2\) - \fracc{\partial_\ttr\p}{\p} |{\overline e}_0^\top|^2.\) - f'' \left|\partial_\ttr^\top\right|^2 \nn \\
&= \(-f'' + \fracc{\cm}{\sm}f'\) \left|\partial_\ttr^\top\right|^2 + f'  \fracc{\cm}{\sm} \(|{\overline e}_0^\top|^2-n\) -  f' \fracc{\partial_\ttr\p}{\p} |{\overline e}_0^\top|^2.
\end{align}
In order that $f'>0$ and the first summand in the last expression be zero, we take $f=\fracc{\cm}{-\mu}$  and the above expression gives
\begin{align}
\(\parcial{}{t} - \Delta\)f(\ttr) &\ge  \cm \(|{\overline e}_0^\top|^2-n\) -  \sm \fracc{\partial_\ttr\p}{\p} |{\overline e}_0^\top|^2   {\ge} \cm \(|{\overline e}_0^\top|^2-n\) -  \eta\ \sm  |{\overline e}_0^\top|^2, \label{evol_f}
\end{align}
where we have used $  \left|\fracc{\partial_\ttr\p}{\p} \right| \le  \left|\fracc{\tN\p}{\p} \right| \le \eta $ for the second inequality.\\
{\bf evolution of $\zeta$}
\begin{align}
\(\parcial{}{t} - \Delta\)\zeta &= e^{\beta n t} \(\parcial{}{t} - \Delta\)f(\ttr) + \beta\ n\  e^{\beta n t}\ f \nn \\
& \us{\eqref{evol_f}}{\ge} e^{\beta n t} \(\cm \(|{\overline e}_0^\top|^2-n\) -  \eta\ \sm  |{\overline e}_0^\top|^2\) + \beta\ n\  e^{\beta n t}\ \fracc{\cm}{-\mu} \label{evol_z}
\end{align}
{\bf evolution of $\phi:= \phi\circ \zeta$}

First, let us observe that $\nabla \phi = \phi' \nabla \zeta$. Then
\begin{align} 
\(\parcial{}{t} - \Delta\)\phi &= \phi' \(\parcial{}{t} - \Delta\) \zeta - \phi'' \left|\nabla\zeta\right|^2 = \phi' \(\parcial{}{t} - \Delta\) \zeta - \phi'' \fracc{|\nabla\phi|^2}{\phi'^2}  \nn \\
&\us{\phi'<0\text{ and} \eqref{evol_z}}{\le} \phi' \(e^{\beta n t} \(\cm \(|{\overline e}_0^\top|^2-n\) -  \eta\ \sm  |{\overline e}_0^\top|^2\) + \beta\ n\  e^{\beta n t}\ \fracc{\cm}{-\mu}\) - \phi'' \fracc{|\nabla\phi|^2}{\phi'^2}  \nn \\
&=\phi' e^{\beta nt} \(\cm |{\overline e}_0^\top|^2-n \cm -  \eta\ \sm  |{\overline e}_0^\top|^2 + \ n\  \ \fracc{\beta}{-\mu} \cm\) - \phi'' \fracc{|\nabla\phi|^2}{\phi'^2}  \lb{evol_ph}
\end{align} 
Now, we choose  
\be
\beta = -\mu  \(1+ \fracc{ \eta}{n}\) , \quad \text{ (notice that $\beta=\nu$ when $ \p=1$)}
\ee
By substitution of this value of $\beta$ in \eqref{evol_ph}, using  $\cm\ge \sqrt{|\mu|} \sm$, we obtain
\begin{align} 
\(\parcial{}{t} - \Delta\) \phi 
& \le  \phi' e^{\beta nt} \cm |{\overline e}_0^\top|^2 - \phi'' \fracc{|\nabla\phi|^2}{\phi'^2}  \le - \phi'' \fracc{|\nabla\phi|^2}{\phi'^2} \label{evol_phi}
\end{align}
and, finaly, the {\bf evolution of $\phi v$}
\begin{align} 
\(\parcial{}{t} - \Delta\)(\phi \ v) &= v \(\parcial{}{t} - \Delta\) \phi +  \phi \(\parcial{}{t} - \Delta\) v - 2 \<\nabla \phi, \nabla v \> \nn \\
&\le - v\ \phi'' \fracc{|\nabla\phi|^2}{\phi'^2} + \phi \(- \frac2{v} |\nabla v|^2     + 2 \<\nabla v, \fracc{\nabla\p}{\p}\> - v  |A|^2 + v \(1-\fracc1{v^2}\) (n-1) \nu \) \nn \\
& \qquad \qquad - 2 \<\nabla \phi, \nabla v\> \nn \\
& = - 2 \a v\  \fracc{|\nabla\phi|^2}{\phi'^2} + 2 \( \<\nabla(\phi v), \fracc{\nabla\p}{\p}\>  - v  \<\nabla \phi, \fracc{\nabla\p}{\p}\> \) \nn \\
& \qquad - \phi  v  |A|^2 + \phi v \(1-\fracc1{v^2}\) (n-1) \nu - 2 \<\nabla( \phi v), \fracc{\nabla v}{v}\> \nn \\
&{\le}  2 \<\nabla(\phi v), \fracc{\nabla\p}{\p} -\fracc{\nabla v}{v}\>  - 2 \a v\  \fracc{|\nabla\phi|^2}{\phi'^2} + 2 \eta v |\nabla \phi|  +  (n-1) \nu  \phi v \label{evol_phiv},
\end{align}
where we have used $\nu \ge 0 \text{ and } \fracc{|\nabla\p|}{\p} \le \fracc{|\tN\p|}{\p} \le \eta, $ for the last inequality. If $\nu < 0$, we can forget the last summand in \eqref{evol_phiv}. To consider both cases in the future we define 
\be\lb{def_enu}
\ev =\left\{\begin{matrix} \nu \text{ if }\nu>0 \\ 0 \text{ if }\ \nu\le 0 \end{matrix}. \right.
\ee

Now, from the obvious inequality $( |\nabla\phi| - 2 \phi'^2 )^2 \ge 0$, we obtain
\begin{align}
2 \eta v |\nabla \phi|  &\le  \fracc{\eta}{2} v \fracc{|\nabla\phi|^2}{\phi'^2} + 2 \eta \phi'^2 v
\end{align}
Then, if we take $\a =  \(1+\fracc{\eta}{4}\)$, we obtain
\begin{align}
\(\parcial{}{t} - \Delta\)(\phi \ v) & < 2  \<\nabla(\phi v), \fracc{\nabla\p}{\p} -\fracc{\nabla v}{v}\>  + 2\ \eta \fracc{\phi'^2}{\phi}\ \phi v + (n-1) \ev\  \phi v
\nn \\
&=  2  \<\nabla(\phi v), \fracc{\nabla\p}{\p} -\fracc{\nabla v}{v}\>  + \( 2\eta^2 + (n-1) \ev \)  \phi v . \label{evol_phiv2} \end{align}

 Let us consider  the sets
$$S_{\r,t} = \{F(x,t)\  ;\  \fracc{\cm}{-\mu}(\r) - \fracc{\cm}{-\mu}(\ttr(F(x,t))) e^{\beta n t} \ge 0\}, \quad \ds \mathbb S_{\r,\tau} = \bigcup_{t\in [0,\tau)} S_{\r,t},  $$ 
(where \lq\lq$)$'' means \lq\lq$]$'' if $\tau<T$ and \lq\lq$[$'' if $\tau=T$).
Observe that $\phi v$ vanishes on the boundary of $S_{\r,t}$ and that 
$$ S_{\r,0} = \{F(x,0;\ \ttr(F(x,0))\le \r\} = \ttr^{-1}(\r) \cap M_0.$$ 
By the maximum principle applied to the function $\phi v$, on $\mathbb S_{\r,T}$, $\phi v(F(x,t))$ is bounded from above by the solution of the equation $y' = (2\eta^2 +(n-1)\ev) y$ with the initial conditions $ \ds y(0)= \sup_{S_{\r,0} } \phi v$, that is,
\begin{multline}
\a \( \fracc{\cm}{-\mu}(\r) - \fracc{\cm}{-\mu}(\ttr(F(x,t))) e^{\beta n t}\)^2 v(F(x,t)) \\
 \le  \sup_{S_{\r,0}}\phi v \  e^{( 2\eta^2 + (n-1) \ev)t}  \\
\le \a \( \fracc{\cm}{-\mu}(\r)\)^2 \sup_{S_{\r,0}} v \   e^{( 2\eta^2 + (n-1) \ev)t}
\label{vbound}
\end{multline}
which, for every $F(x,t)$ in the slice $M_{t},$ taking $\r$ big enough, gives an upper bound of $v(F(x,t))$ (depending on $t$ and $\r$), which finishes the proof of the claim.  
\end{demo}

 \begin{nota}\lb{notagrafo} 
 Let us observe that, unlike in the compact case, the conclusion of Theorem \ref{graphbo2nonco} is just a bound on $v$, with no conclusion about the preservation of the property of being a graph, which will be a consequence of the proof of  the existence theorem in section \ref{STE} and the prolongation theorem \ref{existnonco}. The reason of this is that, when $M_t$ is not compact, it is not clear that $v$ bounded at each point implies $M_t$ is a graph.
 
 However, if $M_t$ is complete and has $|A|$ bounded by an universal constant, it is true that $v$ bounded at each point implies that $M_t$ is a graph over $M$. In fact, let us suppose that $M_t$ is a graph over a proper subset of $M$ and has $|A|$ bounded, let $x\in M$ such that the line $L=\{x\}\times \re$ does not cut $M_t$. The distance $d$ to $L$ is a $C^\infty$ function on $M_t$, and its infimum is the distance $\delta$ between $M_t$ and $L$. Then, by Omori's Lemma (cf \cite{Om}),  there is a sequence of points $p_n\in M_t$ such that $d(p_n)- \delta <\frac1n$ and $|\nabla d|(p_n) < \frac1n$. 
 
  For each $p_n$, let $q_n\in L$ such that $\dist(q_n,p_n) = d(p_n)$. Let $\gamma_n(t)$ be the geodesic from $q_n$ to $p_n$ realizing $d(p_n)$. It is orthogonal to $L$ at $q_n$ and, since the hypersurfaces $u=$constant are totally geodesic (cf. Remark \ref{fh0}), $\gamma_n$ is contained in such a hypersurface and orthogonal to the vector field ${\overline e}_0$. Moreover, it is an integral curve of $\oN  d$ which is then, orthogonal to ${\overline e}_0$. On the other hand  $\nabla d = \oN d-\<\oN d,N\> N$, then $|\nabla d|(p_n) < \frac1n$ implies $\oN d(p_n)$ approach $N(p_n)$ as $n\to\infty$, then, as $n\to\infty$,  $\<{\overline e}_0, N\>_{p_n}$ approach $0$ and $v$ goes to $\infty$. But, taking $\r>\dist(x,x_0)$,  this is in contradiction with \eqref{vbound}. Then, on [0,T[, $M_t$ remains to be a graphic over all $M$.

The property above gives a geometric difference between the concepts of graph in $M\times_\p\re$ and $\re\times_\p M$. In fact,
for graphs in $\re\times_\varphi M$ and $\p$ not constant, it is no longer true that the  hypersurfaces $u=$constant be totally geodesic, then  $\oN d$ is not orthogonal to  ${\overline e}_0$ and the argument fails. It is not dificult, using this idea, to construct a hypersurface in the hyperbolic space $\mathcal H^{n+1}$ which is complete, a geodesic graph on an open set $U$ of $\mathcal H^{n}$ ($U\ne \mathcal H^n$) and with $v$ bounded. For instance, in the Poincar\'e's ball model, take a disc with boundary at the infinite and parallel to the equator.
\end{nota}

\begin{nota}\label{ejapp}
For use in the long time existence theorem, it is convenient to give a more explicit (although less precise) bound for $v$ than that obtained in \eqref{vbound}. For this, first we consider a smaller set $S_{\r,t,\gamma}$ than $S_{\r,t}$, defined, for any positive $\gamma<1$, by 
$$ S_{\r,t,\gamma} = \{F(x,t)\  ;\  \gamma \fracc{\cm}{-\mu}(\r) - \fracc{\cm}{-\mu}(\ttr(F(x,t))) e^{\beta n t} \ge 0\} \subset S_{\r,t}, \quad  \text{ and } S_{\r,\tau,\gamma} = \cup_{t\in [0,\tau)} S_{\r,t,\gamma}.$$ 
From \eqref{vbound} and some obvious inequlities one has, for every $F(x,t) \in \mathbb S_{R,T,\gamma}$,
\begin{align}
\a \( \fracc{\cm}{-\mu}(R) \)^2 (1-\gamma)^2 v(F(x,t)) & \le  \a \( \fracc{\cm}{-\mu}(R)\)^2 \sup_{S_{R,0}} v\  e^{( 2\eta^2 + (n-1) \ev)t},
\end{align}
that is,  for every $F(x,t) \in  S_{R,t,\gamma}$
\begin{align}\lb{vboundex}
  (1-\gamma)^2 v(F(x,t)) \le   \sup_{S_{R,0}} v \   e^{( 2\eta^2 + (n-1) \ev) T}.
\end{align}
Then, if $v$ is bounded on $M_0$, one has $  v(F(x,t)) \le  (1-\gamma)^{-2} \(\sup_{M_0}v\) \  e^{( 2\eta^2 + (n-1) \ev) T}.$ 
Since this formula is true for any $\gamma$ between $0$ and $1$, we have
\be\lb{vboundnoR}
  v(F(x,t))\  \le \  \sup_{M_0}v \  \  e^{( 2\eta^2 + (n-1) \ev) T}.
  \ee
\end{nota}
\begin{notacion}
For the following sections, it will be convenient to introduce the following notation: given any $\r>0$ and $0<\gamma<1$, we define $\r_i$, $i=1,2,3, ... $ by:
\be\label{def_rho}
\gamma\frac{\cm(\r)}{-\mu} = \frac{\cm(\r_1)}{-\mu}, \qquad \gamma\frac{\cm(\r_i)}{-\mu} = \frac{\cm(\r_{i+1})}{-\mu}
\ee
It is simportant to reamrk that $\mathbb S_{\r_i,\tau,\gamma} =\mathbb S_{\r_{i+1},\tau} \subset \mathbb S_{\r_{i},\tau}$ and also that, given any $\r'>0$, $i\in \mathbb N$  and $0<\gamma<1$, there is a $\r$ such that $\r' = \r_i$.
\end{notacion}

\section{Long time existence}\lb{LTE}

First remember the evolution of $|A|^2$, which we take from \cite{CaMi2} (erasing the terms with $\overline H$) because the notation used there is more similar to that of this paper.
\begin{align}\lb{Acu}
 \parcial{|A|^2}{t} &= \Delta |A|^2 - 2|\nabla A|^2 + 2|A|^2\(|A|^2 +\oRic(N,N)\)   \nn
 \\ & \quad -4 \sum_{i,j} \(\oR_{AE_i E_j AE_i E_j}-  \oR_{AE_i E_j E_i AE_j}\)   -2 \sum_{i,j} \(\nablaa_{AE_i} \oR_{N E_j E_i E_j} + \nablaa_{E_j} \oR_{NE_i AE_i E_j}\).
 \end{align}
  
 In our case, using the orthonormal local frame $N, E_1, E_2, ..., E_n$ such that $AE_i = k_i E_i$, we get
\begin{align}\lb{evL2}
 \parcial{|A|^2}{t} &= \Delta |A|^2 - 2|\nabla A|^2 + 2|A|^4 + 2|A|^2 \oRic(N,N)   \nn
 \\ & \quad -4 \sum_{i,j} \(k_i^2 - k_i k_j\) \oR_{E_i E_j E_i E_j}   -2 \sum_{i,j} \ k_i          \(\nablaa_{E_i} \oR_{N E_j E_i E_j} + \nablaa_{E_j} \oR_{NE_i E_i E_j}\) \nn\\ 
 &= \Delta |A|^2 - 2|\nabla A|^2 + 2|A|^4 + 2|A|^2 \oRic(N,N)   \nn
 \\ 
  & \quad -4 \(\sum_{i<j} \(k_i^2 - k_i k_j\) \oR_{E_i E_j E_i E_j}  +\sum_{i<j} \(k_j^2 - k_j k_i\) \oR_{E_j E_i E_j E_i} \) \nn \\
  & \quad -2 \sum_{i} k_i \(\  \tilde{\delta} \oR_N( E_i, E_i) \) \nn \\
  &= \Delta |A|^2 - 2|\nabla A|^2 + 2|A|^4 + 2|A|^2 \oRic(N,N)   \nn
 \\ 
  & \quad -4 \sum_{i<j} \(k_i - k_j\)^2  \oR_{E_i E_j E_i E_j}   -2 
  \< \a  , \tilde{\delta} \oR_N \>.
 \end{align}
where $\a$ is the second fundamental form of $M_t$ and $\tilde{\delta} \oR_N(X,Y):= \sum_j\(\nablaa_X \oR_{N E_j Y E_j} + \nablaa_ {E_j}\oR_{N Y X E_j} \)$.

We cannot deduce any bound for $|A|$ directly  from \eqref{evL2}. Then, as in \cite{EcHu} and many other places  (in fact, working exactly as in \cite{CaMi2}) we shall study the evolution of $\g = (\psi\circ v) |A|^2$, where 
\begin{align}\label{fib}
& \psi(v) := \frac{v^2}{1-\delta v^2}, \\
 \delta&:= \left\{\begin{matrix} & \fracc1{2(\sup_{M_0}{v^2}) e^{2( (n-1)\nu) T}} \qquad & \text{ if  $M$ is compact }\qquad \qquad \qquad \\
& \fracc1{2(\sup_{M_0}{v^2}) e^{2(2\eta^2 + (n-1)\ev) T}}  & \begin{matrix} \text{ if $M$ is complete non-compact } \\ \text{and $v$ is bounded on $M_0$} \end{matrix} \\
& \fracc{(1-\gamma)^4}{2(\sup_{ S_{\r,0}}{v^2}) e^{2(2\eta^2 + (n-1)\ev) T} } &\begin{matrix} \text{ if $M$ is complete non-compact}\\ \text{ and $v$ is not bounded on $M_0$} \end{matrix} \end{matrix}. \right.
\end{align}
From these definitions of $\psi$ and $\delta$ it follows that, for $F(x,t)\in M_t$, $t\in[0,T[$, when $v$ is bounded on $M_0$ and for $F(x,t)\in \mathbb S_{\r,T,\gamma}$ in the other case, 

\be\lb{pordeff}
 \fracc1{1-\delta} \le \psi \le \fracc1\delta \text{ , } \(- \frac{2}{v} + \delta  \) \fracc1{ \psi'}  - \frac{\psi''}{{(\psi')}^2} +\frac3{2 \psi}   < 0 \text{ 
, }
 - v  \fracc{\psi'}{\psi^2}  + 2 \fracc1{\psi} = -2 \delta <0,
 \ee
 \be\lb{pordeff2}
 \delta \le \fracc12 \text{ , } \fracc{\sqrt 2}{(1- \delta) \sqrt\delta}  \le \fracc{\psi'}{\psi} \le 4  \text{ 
, }
 \fracc{\sqrt 2}{(1- \delta) \delta^{3/2}}  \le \fracc1\delta \fracc{\psi'}{\psi} \le \fracc4\delta \text{     and } 0 \le  \fracc{v^2-1}{v} \fracc{\psi'}{\psi} \le 4 {(1 - 2\delta)}.
\ee

\begin{lema}\lb{l_evolg} For every $F(x,t)\in F(M\times[0,T[)$ if $v$ is universally  bounded on $M_0$, and for every  $F(x,t) \in \ds \mathbb S_{\r_1,T}$ in other case,  one has:
\begin{align}
&\(\parcial{}{t} -  \Delta\) \g
\le   - \frac1{\psi} \<\nabla \g, \nabla\psi\>     -2 \delta \g^2 + \frak B \ \g + \frak C \  \sqrt\g, \quad  \text{ where }\nn 
\\  & \frak B = \fracc1\delta \fracc{\psi'}{\psi}    \(\fracc{|\tN\p|}{\p}\)^2    -   \fracc{v^2-1}{v}  \(\frac{\tD \p}{\p}  + \oRic_{\t1\t1} +  \fracc{|\tN \p|^2}{\p^2}\)  \ \fracc{\psi'}{\psi} + 2  \oRic(N,N)   + 8   |\overline{Scal}-\oRic(N,N)| ,   \nn \\
&    \frak C = 2 \sqrt{\psi} \  |\tilde{\delta} \oR_N |  .
\label{evol_gb}
 \end{align}
\end{lema}
\begin{demo}
The evolution of $\g$ is given by those of  $\psi$ and $|A|^2$ according to the formula:
\begin{align} \lb{fegb}
&\(\parcial{}{t} -  \Delta\) \g = |A|^2 \(\parcial{}{t} - \Delta\) \psi +  \psi \(\parcial{}{t} - \Delta\) |A|^2 - 2 \< \nabla \psi, \nabla |A|^2\> \nn \\
& \qquad = |A|^2 \psi' \(\parcial{}{t} - \Delta\) v - |A|^2 \psi'' \  |\nabla v|^2 + 
 \psi     \(\parcial{}{t} - \Delta\) |A|^2 - 2 \< \nabla \psi, \nabla |A|^2 \>  .
\end{align}
For the last summand in \eqref{fegb}, we use the following inequality (cf. \cite{EcHu}, page 555, but be aware that in \cite{EcHu} is used $\p(v^2)$ instead of $\psi(v)$.)
\be\lb{EcHub}
- 2 \< \nabla\psi , \nabla |A|^2 \> \le - \frac1{\psi} \<\nabla \g, \nabla\psi\> + 2 \psi |\nabla|A||^2 + \frac3{2 \psi} |A|^2 |\nabla\psi|^2
\ee
and Kato's inequality
\be\label{Katb}
| \nabla |A| |^2 \le |\nabla A|^2. (\text{ equivalently } \left|\nabla |A|^2\right|^2 \le 4 |A|^2 |\nabla A|^2)
\ee
From \eqref{fegb}, \eqref{EcHub} and \eqref{Katb}, 
\begin{align}
\(\parcial{}{t} -  \Delta\) \g &\le |A|^2 \psi' \( - \frac2{v} |\nabla v|^2  + \fracc2{\p} \<\nabla v, \nabla\p\>  - v  |A|^2 \right. \nn \\
& \qquad \qquad  \quad \left.   - v \(1-\fracc1{v^2}\) \(\frac{\tD \p}{\p}  + \tRic_{\t1\t1} +  \fracc{|\tN \p|^2}{\p^2}-\fracc{\tN^2\p}\p  \(\t1,\t1\)\) \)
 \nn \\
&- |A|^2 \psi''\ |\nabla v|^2  \nn \\
& +  \psi \( - 2|\nabla A|^2 + 2|A|^4 + 2|A|^2 \oRic(N,N)    -4 \sum_{i<j} \(k_i - k_j\)^2  \oR_{E_i E_j E_i E_j}   -2 
  \< \a  , \tilde{\delta} \oR_N \>\)  \nn \\ 
  &     - \frac1{\psi} \<\nabla \g, \nabla\psi\> + 2\ \psi \  |\nabla A|^2 + \frac3{2 \psi} |A|^2 |\nabla \psi|^2.
\end{align}
Since $\psi'>0$,  using Young's inequality $x y \le \varepsilon x^2 + \frac1{4 \varepsilon} y^2$ with $x= 2 \fracc{|\oN\p|}{\p}$, $y = |\nabla v|$ and 
$\varepsilon:= \fracc1{4 \delta} $, we conclude that
\begin{align}
\(\parcial{}{t} -  \Delta\) \g
\le & \ |A|^2 \psi' \( - \frac2{v} |\nabla v|^2 +  \fracc1\delta \(\fracc{|\oN\p|}{\p}\)^2 + \delta  |\nabla v|^2  - v  |A|^2 \right. \nn \\
& \qquad \qquad  \quad \left.   - v \(1-\fracc1{v^2}\) \(\frac{\tD \p}{\p}  + \tRic_{\t1\t1} +  \fracc{|\tN \p|^2}{\p^2}-\fracc{\tN^2\p}\p  \(\t1,\t1\)\) \)
 \nn \\
&- |A|^2 \psi''\ |\nabla v|^2  \nn \\
& +     2 \ \psi \ |A|^4 + 2\ \psi\ |A|^2 \oRic(N,N)    - 4 \ \psi\ \sum_{i<j} \(k_i - k_j\)^2  \oR_{E_i E_j E_i E_j}   + 2 \ \psi\  
  |A|  |\tilde{\delta} \oR_N|  \nn \\ 
  &     - \frac1{\psi} \<\nabla \g, \nabla\psi\>  + \frac3{2 \psi} |A|^2 |\nabla \psi|^2.
  \lb{evolg1b}
\end{align}

Next, let us bound and/or rearrange the different terms in \eqref{evolg1b}.
\begin{align}
 \(- \frac{2}{v}+\delta\) & |A|^2 \psi'  |\nabla v|^2  - |A|^2 \psi'' |\nabla v|^2 + \frac3{2 \psi} |A|^2 |\nabla \psi|^2 \nn \\
 &= \(\(- \frac{2}{v} +\delta  \) \fracc1{ \psi'}  - \frac{\psi''}{{(\psi')}^2} +\frac3{2 \psi}\) |A|^2 |\nabla \psi|^2 \le 0 \text{ by \eqref{pordeff} }\lb{bdphi2}
\end{align}
\begin{align}
& \ |A|^2 \psi'  \fracc1\delta  \(\fracc{|\tN\p|}{\p}\)^2    -  \(1- \fracc1{v^2}\)  \(\frac{\tD \p}{\p}  + \tRic_{\t1\t1} +  \fracc{|\tN \p|^2}{\p^2}-\fracc{\tN^2\p}\p  \(\t1,\t1\)\) v\ |A|^2 \psi'  \nn \\
&-   |A|^4  v  \psi'    + 2 \psi |A|^4 + 2\ \psi |A|^2 \oRic(N,N)    - 4 \psi  \sum_{i<j} \(k_i - k_j\)^2  \oR_{E_i E_j E_i E_j}   + 2 \psi
  | A | \  |\tilde{\delta} \oR_N |    \nn \\
& \le \(\fracc{\psi'}{\psi}  \  \fracc1\delta \(\fracc{|\tN\p|}{\p}\)^2    -  \(1- \fracc1{v^2}\)   \(\frac{\tD \p}{\p}  + \tRic_{\t1\t1} +  \fracc{|\tN \p|^2}{\p^2}-\fracc{\tN^2\p}\p  \(\t1,\t1\)\)  v\ \fracc{\psi'}{\psi}\) \g   \nn \\
& + \(-     v  \fracc{\psi'}{\psi^2}  + 2 \fracc1{\psi}\) \g^2 \    + 2  \oRic(N,N)  \g  + 8   |\overline{Scal}-\oRic(N,N)|  \g + 2 \sqrt{\psi} \  |\tilde{\delta} \oR_N |  \sqrt{\g}  .
\end{align}
because $- 4 \psi  \sum_{i<j} \(k_i - k_j\)^2  \oR_{E_i E_j E_i E_j}  \le 4 \psi  4 |A|^2 \sum_{i<j}   \oR_{E_i E_j E_i E_j}  \le 2 \psi 4 |A|^2 |\overline{Scal}-\oRic(N,N)|$. And the formula of the theorem follows putting all this together and using the relation between $\tR$ and $\oR$ computed in section \ref{curvoM}.
\end{demo}

\begin{teor}\label{existco} 
Let $M$ be compact. If $\oM = M \times_{\p}\re$ and $M_0$  is a $C^\infty$ graph over $M$, then   \eqref{vpmf}  has a solution whith initial condition $M_0$, defined on a maximal time interval $[0,\infty[$ and which is a graph for all $t$. 
More, when $\nu \le 0$,  $|A|$ is bounded on $[0,\infty[$ with a bound which does not depend on $t$.
\end{teor}
\begin{demo} 
From \eqref{pordeff} and \eqref{pordeff2} we obtain the following inequalities and define the constants $K$ and $C$ by
\begin{align} 
\frak B &\le   \fracc4\delta \eta^2    +  4  (1-2\delta) \  \ev   + (2 + 8 n) \sup_\oM |\oRic| =: K  \lb{K}
\end{align}
\begin{align}
\frak C \le    \fracc4{\sqrt{\delta}}\  \sup_\oM |\oN \oR |  =: C . \lb{C}
\end{align}
From these definitions of $K$ and $C$ and Lemma 4 it follows
\begin{align}
\(\parcial{}{t} -  \Delta\) \g
\le &  - \frac1{\psi} \<\nabla \g, \nabla\psi\>     -2 \delta \g^2 + K \g + C \sqrt{\g} \label{evol_g}
 .\end{align}

We remark that $K$ and $C$  depend on $T$ only through $\delta$, $\delta$ depends on $\sup v$ and, according to Theorem \ref{graphbo2comp}, when $M$ is compact and $\nu \le 0$, $\sup v$ does not increase  with time.

As we  remarked in section \ref{evoluT}, when $M$ is compact,  there is a  solution of \eqref{vpmf} for a maximal interval $[0,T[$. Let $t_1\in[0,T[$. Let us suppose that $(x_0,t_0)$  is the point where $\g$ attains its maximum for $t\le t_1$ and $0< t_0 <t_1$.

If $\g_0:=\g(x_0,t_0) >1$, then  it has to satisfy $0 \le - 2 \ \delta(t_0) \g_0^2  + (K + C)(t_0)\ \g_0  \le - 2 \ \delta \g_0^2  + (K + C)\ \g_0$, so 
\be\lb{ubg}
\g_0 \le \max\{\max_{M_0} \g,  \fracc{K+C}{2\delta}, 1\}.
\ee
Since $\psi$  and $\g$ are bounded, $ |A|^2$ is bounded on $[0,T[$ with a bound which does not depends  on $T$ when $\nu \le 0$, because in this case neither $K$, $C$ and $\delta$ depend on $t$.

Once we achieve the upper bound for  $|A|^2$, it  follows, like in  \cite{Hu84} and \cite{Hu86}, that $|\nabla^jA|$ is bounded for every $j\ge 1$. If $T<\infty$, these bounds imply (cf. \cite{Hu84} pages 257, ff.) that   $X_t$ converges (as $t\to T$, in the $C^\infty$-topology) to a unique smooth  limit $X_T$. Now we can apply  the short time existence theorem to continue the solution after $T$,  contradicting the maximality of $[0,T[$.
\end{demo}

The next is a theorem modulo the existence theorem that will be proved in section \ref{STE}. For this reason, in its hypotheses appears the existence of a short time existence theorem.

\begin{teor}\label{existnonco} 
Let $M$ be complete non compact with a pole $x_0$.   Let $\oM = M \times_{\p}\re$. Let us suppose that for every graph $\mathcal M_0$ over $M$ there is a solution of \eqref{vpmf} with initial condition $\mathcal M_0$ that is a graph and defined on some time interval. Let $M_t$ be a solution of  \eqref{vpmf}  defined on a maximal time interval $[0,T[$. If  $M_0$ is a $C^\infty$ graph over $M$, then   $T=\infty$. 
\end{teor}
\begin{demo} Let $\eta$, $\nu$ and $\ev$ be as defined by \eqref{def_mueta}, \eqref{def_nu} and \eqref{def_enu}. For the terms $\frak B$ and  $\frak C$ in Lemma \ref{l_evolg} we have the same bounds that in  Theorem \ref{existco}, so we define $K$ and $C$ again by \eqref{K} and \eqref{C} respectively. Now we shall use the functions $\phi$ and $\zeta$ defined in Theorem \ref{graphbo2nonco} and  compute the evolution of $\phi \g$ for points in $S_{\r_1,T}$.
\begin{align} 
\(\parcial{}{t} - \Delta\)(\phi \ \g) &= \g \(\parcial{}{t} - \Delta\) \phi +  \phi \(\parcial{}{t} - \Delta\) \g - 2 \<\nabla \phi, \nabla \g \> \nn \\
& \le - 2 \a \g \fracc{|\nabla\phi|^2}{\phi'^2} + \phi \(-\fracc1{\psi} \<\nabla\g.\nabla\psi\> - 2 \delta \g^2 + K \g + C \sqrt{\g}\) - 2 \<\nabla\phi, \nabla \g\> \nn\\
& = - 2 \a \g \fracc{|\nabla\phi|^2}{\phi'^2}   - \<\nabla(\phi\g) - \g \nabla \phi, \fracc{\nabla\psi}{\psi}\> - 2 \frac{\delta}{\phi} (\phi \g)^2 + K \phi\g \nn \\
& \qquad \qquad + C \sqrt{\phi}\sqrt{\phi \g} - 2 \<\frac{\nabla\phi}{\phi}, \nabla (\phi \g) - \g \nabla\phi\> \lb{phg0}
\end{align}
but $\phi'^2 = 4 \a \phi$, then $- 2 \a \g \fracc{|\nabla\phi|^2}{\phi'^2} + 2 \g \fracc{|\nabla\phi|^2}{\phi} = 6 \a \fracc{|\nabla\phi|^2}{\phi'^2} \g= 6 \a |\nabla\zeta|^2 \g= 6 \a \g \sm^2 |\partial_{\ttr}^\top|^2e^{2 \beta n t}$. By substitution in the evolution of $\phi \g$, we obtain
\begin{align} 
\(\parcial{}{t} - \Delta\)(\phi \ \g) & \le 6 \a \g \sm^2 e^{2 \beta n t}   - \<\nabla(\phi\g) - \g \nabla \phi, \fracc{\nabla\psi}{\psi}\> - 2 {\delta} \  \phi\   \g^2 + K \phi \g \nn  \\
& \qquad \qquad + C \sqrt{\phi}\sqrt{\phi \g} - 2 \<\frac{\nabla\phi}{\phi}, \nabla (\phi \g) \> \nn \\
&= 6 \a \g \sm^2 e^{2 \beta n t}   - \<\nabla(\phi\g),   \fracc{\nabla\psi}{\psi} + 2 \frac{\nabla\phi}{\phi}\> + \<\g \nabla \phi, ,   \fracc{\nabla\psi}{\psi} \>\nn  \\
& \qquad \qquad - 2 \delta \ \phi\ \g^2 + K \phi \g  + C \sqrt{\phi}\sqrt{\phi \g}  .\lb{phg1}\end{align}
Using $\nabla \psi = \psi' \nabla v$ and the expression \eqref{nablav} of $\nabla v$,
\begin{align}
\fracc{\nabla \psi}{\psi} = \frac{1-\delta^2 v^4}{v^2} \nabla v = (1-\delta^2 v^4)\(\<\fracc{\tN\p}{\p},N\> {\overline e}_0^\top + A{\overline e}_0^\top\) , \text{ then }
\end{align}
\begin{align}
\left|\fracc{\nabla \psi}{\psi} \right| \le \eta + |A| = \eta + \fracc{\sqrt{\g}}{\sqrt{\psi}} .\end{align}
On the other hand
\begin{align}
\sm^2 e^{2\beta n t} \le \(\cm e^{\beta n t}\)^2 \le  \cm(\r_1)^2 \text{ on } \mathbb S_{\r_1,T} .
\end{align}
Pluging these inequalitues in \eqref{phg1},
\begin{align} 
\(\parcial{}{t} - \Delta\)(\phi \ \g) & \le  6 \a  \cm(\r_1)^2 \g   - \<\nabla(\phi\g),   \fracc{\nabla\psi}{\psi} + 2 \frac{\nabla\phi}{\phi}\> +  |\nabla \phi| \(\eta +  \fracc{\sqrt{\g}}{\sqrt{\psi}}\) \g\nn  \\
& \qquad \qquad - 2 \delta \ \phi\ \g^2 + K \phi \g  + C \sqrt{\phi}\sqrt{\phi \g}  .\nn \\
& = - \<\nabla(\phi\g),   \fracc{\nabla\psi}{\psi} + 2 \frac{\nabla\phi}{\phi}\> - 2 \delta \ \phi\ \g^2 \nn \\
&\qquad + \fracc{\nabla \phi}{\sqrt{\psi}} \g \sqrt{\g} +  \( 6 \a  \cm(\r_1)^2 + \eta |\nabla\phi|  + K \phi\) \g + C \sqrt{\phi}\sqrt{\phi \g} \lb{phg2}  \end{align}
At a point where a maximum for $\phi\g$ is attained it follows from the above inequality that, if this happens when $t\ne 0$,
\be
2 \delta \ \phi\ \g^2 \le \fracc{\nabla \phi}{\sqrt{\psi}} \g \sqrt{\g} +  \( 6 \a  \cm(\r_1)^2 + \eta |\nabla\phi|  + K \phi\) \g + C \sqrt{\phi}\sqrt{\phi \g} , \nn
\ee
and, multiplying by $\fracc{\sqrt\phi}{\sqrt\g }$,
\be
2 \delta \(\sqrt{ \phi \g }\)^3 \le \fracc{\nabla \phi}{\sqrt{\psi}\sqrt\phi}  \phi\g+   \(6 \a  \cm(\r_1)^2 + \eta |\nabla\phi|   + K \phi \) \sqrt{\phi \g}+ C  \phi\sqrt{\phi}. \lb{phg5}
\ee
 having into account that
\begin{align}
\nabla \phi = - 2  \sqrt{\a \phi}\  \sm(\ttr) \partial_{\ttr}^\top, \quad \phi \le \a \fracc{\cm(\r_1)^2}{\mu^2}, \quad  \quad \fracc1{\sqrt\psi} \le \sqrt{1-\delta}, \lb{insqph}
\end{align}
we have 
\begin{multline}
2 \delta \( \sqrt{\phi\ \g}\)^3 \le 2 \sqrt{\a} \sm(\r_1) \sqrt{1-\delta}   \( \sqrt{\phi\ \g}\)^2  
\\ +  \( \( 6 +\fracc{K}{\mu^2} \) \a\ \cm(\r_1)^2 + \eta 2 \a \fracc{\cm(\r_1)}{-\mu} \sm(\r_1) \)   \sqrt{\phi\ \g} + C \a^{3/2} \fracc{\cm(\r_1)^3}{-\mu^3} .  \nn
\end{multline}
Having into account that $\sm(\r_1)\le \cm(\r_1)$ and dividing by $\cm(\r_1)^3$, we get
\begin{multline}
2 \delta \( \fracc{\sqrt{\phi\ \g}}{\cm(\r_1)}\)^3 - 2 \sqrt{\a}  \sqrt{1-\delta}   \( \fracc{\sqrt{\phi\ \g}}{\cm(\r_1)}\)^2  
-  \( \( 6 +\fracc{K}{\mu^2} \) \a\  + \eta 2 \a \fracc{1}{-\mu}  \)   \fracc{\sqrt{\phi\ \g}}{\cm(\r_1)}  - C \a^{3/2} \fracc{1}{-\mu^3} \le 0.  \lb{phg4} 
\end{multline}
As the coefficient of $ \( \fracc{\sqrt{\phi\ \g}}{\cm(\r_1)}\)$ is positive, the above inequality implies  that 
\be \( \fracc{\sqrt{\phi\ \g}}{\cm(\r_1)}\) \le D=\max\{\mathcal D,\max_{S_{\r_1,0}}\sqrt{\phi\g}\}, \lb{phg3}
\ee 
where $\mathcal D$ is the biggest solution of the third order polynomial 
equation given by \eqref{phg4} when we change the inequality by an equality. $\mathcal D$ depends on the coefficients of the equation, then it only depends on the geometry of $\oM$ and $\delta$.
From the above inequality we have that,  every $F(x,t) \in \mathbb S_{\r_1,T,\gamma}$  one has
\be
    \(\fracc{\sqrt{ \a} (1-\gamma)}{-\mu}\) \sqrt{\g}  \le D, \nn
   \ee
then, remembering that $\sqrt\g = \sqrt\psi \ |A|$ and the last inequality \eqref{insqph},
  \be
  |A| \le \fracc{- \mu \sqrt{(1-\delta)}  D}{\sqrt\a (1-\gamma)}, \label{oo}
  \ee 
   This shows that, on $\mathbb S_{\r_2,T} =  \mathbb S_{\r_1,T,\gamma}$, $|A|$ is bounded by a bound $C(0,\r, \gamma, T)$ depending   $T$, $\gamma$ and  $\rho$ (only through $\delta$) and on  $M_0$. 
      
   From here, a variant of  the  procedure used by Ecker, Huisken and Unterberger works. Let us give some details of  these computations.
Let us suppose, by induction, that, on $\mathbb S_{\r_{2+k},T}$,
 \be\lb{bo_dtAk}
 |\nabla^k A| \le C(k,\r,\gamma,T) \quad \text{ for } \quad  k=0, ..., m-1\ee
where $C(k,\r,\gamma,T)$ depends on $\r$, $\gamma$ and $k$ (through $\rho_{2+k}$), the bounds of $|\nabla^iA|$ for $0\le i \le k-1$ and the geometry of $\oM$ (the bounds on $|\oN^j\oR|$, $j=0, ..., k$).
 
  Following the same procedure of \cite{Hu84} and \cite{Hu86}, we can find a constant $D_1$ depending on $m$, $n$, the geometry of $\oM$ (that is, the bounds on $|\oN^k\oR|$, $k=0, ..., m$) and the bounds of $|\nabla^kA|$, $k=0, ..., m-1$ such that 
\begin{align}
\parcial{}{t} |\nabla^{m} A |^2 &\leq \Delta |\nabla^{m} A |^2   - 2 |\nabla^{m+1} A|^2+  D_1 (|\nabla^{m} A |^2 + 1) \lb{dtAm1}
\end{align}

Next, we define  
\be \lb{deff}
f:= |\nabla^{m} A |^2 + \xi |\nabla^{m -1} A |^2,
\ee
 being $\xi$ a constant to be specified later. For the time derivative of $f$, \eqref{dtAm1} yields
\be \lb{derf_aux}
\parcial{f}{t} \leq \Delta |\nabla^{m} A |^2  + D_1 \(|\nabla^{m} A |^2 + 1\) + \xi \parcial{}{t} |\nabla^{m -1} A |^2
\ee

Observe that the last addend on the right hand side of \eqref{derf_aux} can be estimated using again \eqref{dtAm1}. In fact,
\begin{align*}
\parcial{|\nabla^{m -1} A |^2}{t} & = \Delta |\nabla^{m -1} A |^2 -2 |\nabla^{m} A |^2 + D_2,
\end{align*}
being 
$D_2$ a constant depending on $C(m-1,\r,\gamma,T)$.

Substituting this in \eqref{derf_aux}, we get 
\begin{align*}
\parcial{f}{t}  & \leq \Delta f + (D_1 - 2 \xi) |\nabla^m A|^2 + D_1 + \xi D_2.
\end{align*} 
If we choose $\xi \geq D_1$ and $D_3:= D_1 + \xi D_2$, then
\begin{align*}
\parcial{f}{t} & \leq \Delta f -  \xi |\nabla^{m} A |^2 + D_3 \us{\eqref{deff}}{=}   \Delta f  -  \xi f + \xi^2 |\nabla^{m - 1} A|^2 + D_3 
\\ & \leq \Delta f -  \xi f + \xi^2 C_{m -1} + D_3.
\end{align*}
Considering the same function $\phi$ as before, the same computational rules and the evolution of $\phi$, we get
\begin{equation} 
\(\parcial{}{t} - \Delta\)(\phi \ f) \le  - 2 \a f \fracc{|\nabla\phi|^2}{\phi'^2}   - \xi\  \phi f +( \xi^2 C_{m -1} + D_3) \phi   - 2 \<\frac{\nabla\phi}{\phi}, \nabla (\phi f) - f \nabla\phi\> \lb{phg0f}
\end{equation}
The same computations than in the proof of Theorem \ref{existnonco} give now 
\begin{align} 
\(\parcial{}{t} - \Delta\)(\phi \ f) & \le  6 \a  \cm(\r_{m+1})^2 f   - \<\nabla(\phi f),   2 \frac{\nabla\phi}{\phi}\>   - \xi\  \phi f +( \xi^2 C_{m -1} + D_3) \phi     \lb{phg2f}  \end{align}
At a point of $\mathbb S_{\r_{m+1},T}$  where a maximum for $\phi f$ is attained it follows from the above inequality that
\be
( \xi\  \phi   -    6 \a  \cm(\r_{m+1})^2)  \phi f    \le ( \xi^2 C_{m -1} + D_3) \phi^2 \le ( \xi^2 C_{m -1} + D_3) \fracc{\a^2}{\mu^4} \cm(\r_{m+1})^4,
\ee
then, at any point in $\mathbb S_{\r_{2+m},T}=\mathbb S_{\r_{m+1},T,\gamma}$, since, on these points,  $\phi\ge \fracc{\a}{\mu^2} \cm(\r_{m+1})^2 (1-\gamma)^2$,
\begin{align*}
( \xi\ \fracc{\a}{\mu^2}  \cm(\r_{m+1})^2 (1-\gamma)^2    -    6 \a  \cm(\r_{m+1})^2)  \phi f   & \le ( \xi^2 C_{m -1} + D_3) \fracc{\a^2}{\mu^4} \cm(\r_{m+1})^4,
\end{align*}
If we take $\xi =\max \{D_1,  \fracc{7  \mu^2}{(1-\gamma)^2}\}$, we have, for points in  $\mathbb S_{\r_1,T,\gamma}$,
\be
\a  \cm(\r_{m+1})^2 \fracc{\a}{\mu^2} \cm(\r_{m+1})^2 (1-\gamma)^2 f  \le  ( \xi^2 C_{m -1} + D_3) \fracc{\a^2}{\mu^4} \cm(\r_{m+1})^4,
\ee
that is
\be\lb{oom}
 f  \le ( 1-\gamma)^{-2} ( \xi^2 C_{m -1} + D_3) \fracc{1}{\mu^2}=: C(m,\r,\gamma,T) ,
\ee

Finally,  by \eqref{deff}, $|\nabla^{m} A|^2 \leq f \le C(m,\r,\gamma,T)$ on $\mathbb S_{\r_{2+m},T}$.

Now the final argument cannot be exactly like in the compact case because there, to bound $\dist(F(p,t),x_0)$, one uses integration along $t$, but here it could happen that $F(p,t) \in S_{\r_{2+m},T}$ but $F(p,s) \notin S_{\r_{2+m},T}$ for $s<t$, then we do not know anything about the bounds of $|A|$ in $F(p,t)$.  To avoid this problem, we have to use the parametrization $\overline F(\cdot,t)$ of $M_t$ and the equation \eqref{mcf} or its equivalent \eqref{mcft}.  From the bounds on $v$ obtained in the section \ref{graph} and \eqref{nuv} we get that the gradient of $u$ is bounded on $S_{\r_1,T}$. From this, formula \eqref{Aeiej} and the above  bounds on $|\nabla ^mA|$ we get that  the higher order derivatives of $u$ are bounded. And also from the bounds on $v$ and $|A|$ and formula \eqref{mcfg} it follows that $|u(t)| \le |u(0)|+ \int_0^T \frac{v}{\p} n |A|$ is bounded on $S_{\r_2,T}$. Once we have these bounds, we have, on each compact $\left\{x\in M_t;  \fracc{\cm}{-\mu}(\r_{2+m}) - \fracc{\cm}{-\mu}(\ttr(x)) e^{\beta n t} \ge 0\right\}$ a well defined limit of $M_t$ when $t\to T$, which allows to continue the flow after $T$. Then $T = \infty$.
\end{demo} 

\section{Existence of solution when $M$ is non-compact}\lb{STE}

When $M$ is complete non-compact and posses a pole $x_0$, the existence of solution follows from the parabolic theory using the estimates of section \ref{graph} and \ref{LTE} in the following way. 

Given $\r_0>0$, let us define $\r=\r(\r_0)$ by the expression
\be\lb{R'R0}
 \frac{\cm(\r_{2+m})}{-\mu}  e^{-\beta n T} = \frac{\cm(\r_0)}{-\mu}.\ee
For any $\r' \ge \r$, by the theory of parabolic equations, there is a unique solution $u_{\r'}$ of \eqref{mcft} on $\widehat B(x_0,\r') \times [0,T]$ satisfying the conditions:
\begin{align}
u_{\r'}(\cdot,0)= u_0|_{\widehat B(x_0,{\r'})} \quad \text{ and } \quad 
u_{\r'}(x,t) = u_0(x) \text{ for } (x,t) \in \partial \widehat B(x_0,{\r'}) \times [0,T].
\end{align}
From \eqref{R'R0}  one has that if $(x,t)\in \widehat B(x_0,\r_0) \times [0,T]$, then $(x,u_{\r'}(x,t)) \in S_{\r_{2+m},T}$ and, from the estimates \eqref{vboundex}, \eqref{oo} and \eqref{oom} we get that there are constants $c(m,\r,\gamma,T)$ such that $|\tN^m u_R| \le c(m,\r,\gamma,T)$ and, since $\r$ is determined from $\r_0$  and $\gamma$ by \eqref{R'R0}, we can say that $|\tN^m u_R|$ are bounded by constants $C(m,\r_0,\gamma,T)$.
From these bounds one gets also an estimate for $|u_\r'|$ in the same way that was done at the end of the proof of Theorem \ref{existnonco}.

Now, observe that, thanks to \eqref{mcfg} and \eqref{bo_dtAk}, $\ds \left|\parcial{u}{t}\right|$ is bounded  by a constant depending only on $m$, $\r_0$ $\gamma$ and $T$. 
From \eqref{mcfg}, $\ds\left|\parcial{^2u}{t^2}\right|$  is bounded if $v$, $|H|$, $\left|\parcial{v}{t}\right|$ and $\left|\parcial{H}{t}\right|$ are bounded. We proved before that the first two quantities are bounded. The bounds on the other two follow from \eqref{bo_var_v}, \eqref{var_HtT}, Proposition \ref{propax}, formula \eqref{nablav} and the fact (proven before) that $v$, $|A|$, $|\nabla A|$ , $|\nabla^2 A|$, $|\tN u|$ and $|\tN^2u|$  are bounded.

As a consequence, for every $\r_0>0$, the family of $C^\infty$ functions $\{u_{\r'}\}_{\r'\ge \r}$ converges to a smooth function $\widehat u_{\r_0}$ on $\widehat B(x_0,\r_0)$  which is at least $C^1$ on $t$ and a solution  of \eqref{mcft} on $\widehat B(x_0,\r_0) \times [0,T]$ (then, by parabolic theory, it is also $C^\infty$ on $t$). 
Given a sequence $\r_0^1 < \r_0^2 < \cdots \to \infty$, for $j>i$ the families  $\{u_{\r'}\}_{\r'\ge \r(\r_0^i)}$ and $\{u_{\r'}\}_{\r'\ge \r(\r_0^j)}$,  coincide for $\r'\ge \r(\r_0^j)$, then  their  limits  $\widehat u_{\r^i}$ satisfy the property that $\widehat u_{\r^j}|_{\widehat B(x_0,\r^i)} = \widehat u_{\r^i}$, then they define a smooth function $\widehat u$ on $M$ which is the $C^\infty$ limit on the compacts of the family $u_{\r'}$ when $\r'\to\infty$, and it is a solution of \eqref{mcft}. Then, joining these arguments with those of the proof of Theorem \ref{existnonco} we have

\begin{teor} \label{stenonco} 
Let $M$ be complete non compact with a pole $x_0$.  If $\oM = M \times_{\p}\re$ and $M_0$ is  a $C^\infty$ graph over $M$, then there is a  solution $M_t$ of  \eqref{vpmf}   with initial condition $M_0$ which is a graph over $M$ and is defined on $[0,\infty[$. 
\end{teor}

\section{Existence of solution for a Lipschitz initial condition}\lb{STEL}

The existence of solution of \eqref{vpmf} when $M_0$ is a graph over $M$ given by a Lipschitz continuous function  follows by approximating $M_0$ by a sequence of smooth graphs $M_n$ and applying the existence theorems \ref{existco} and \ref{stenonco} to these approximations. But, in order to show that the solutions of \eqref{vpmf} with initial conditions $M_n$ converge to a smooth solution, we need to get bounds of $|\nabla^mA|$ that do not depend on the bounds on the initial condition, because these could go to $\infty$ as $n\to\infty$ because the limit $M_0$ of $M_n$ is only Lipschitz. In this section we shall obtain these estimates.

\begin{lema}\lb{boundNAmL} Given a smooth solution of \eqref{vpmf} defined on $[0,T[$, for every $m=0,1,2,...$ there is a constant $\a_m$ such that
 \be \lb{evol_NAmt}
 |\nabla^m A|^2 \le \alpha_m \(\fracc{t+1}{t}\)^{m+1}
\ee
on $M_t$ if $M$ is compact or in $\mathbb S_{\r_{2+m},T}$ if $M$ is non-compact. Moreover, the constant $\a_m$ depends  on $m$, the geometry of $\oM$ and $T$ in the first case, and also on  $\r$ and $\gamma$ in the second.
\end{lema}
\begin{demo}
First, let us consider the case $M$ compact. For obtaining the bound when $m=0$ we start using the inequality \eqref{evol_g} to obtain the following inequation for the evolution of $t \frak g$
\be\lb{evol_tg}
\(\parcial{}{t} - \Delta\) (t \frak g) \le -\fracc1{\psi} \<\nabla(t \frak g), \nabla \psi\> - 2 \ \delta \ \frak g^2\ t + K\ \frak g \  t  + C\ \sqrt{\frak g} \ t + \frak g.
\ee
Given $t \in ]0,T[$, let $t_0$ be the time when $t \frak g$ attains it maximum value $t_0 \frak g_0$ in $[0,t] \times M$ (then $\frak g_0 = \max_{x\in M_{t_0}} \frak g(x)$). By the maximum principle, from \eqref{evol_tg} we get 
\be\lb{mp_tg}
 2 \ \delta \ \frak g_0^2\ t_0 \le   (1 + K\ t_0) \frak g_0  + C\ \sqrt{\frak g_0} \ t_0.
\ee

If $\frak g_0\le 1$, then, by definition of maximum, $t\ \frak g \le t_0 \frak g_0 \le t_0 \le t$, then $\frak g\le 1$ on $M_t$, and, by the definition of $\frak g$ and \eqref{pordeff},  $|A|^2 = \fracc1{\psi} \frak g \le (1-\delta) \frak g  \le 1-\delta  \le \fracc{1+t}{t}$.

If $\frak g_0 \ge 1$, then $t_0 \sqrt{\frak g_0} \le t_0 \frak g_0$ and it follow from \eqref{mp_tg} and the definition of maximum that
\be 
2 \delta \frak g t \le 2 \delta \frak g_0 t_0 \le (1+ (K+C) t_0)  \le (1+ (K+C) t), 
\ee
then, using again  \eqref{pordeff},
\be
|A|^2 = \fracc1{\psi} \frak g \le (1-\delta) \fracc{1+ (K+C) t}{2 \delta t} \le \a_0 \fracc{1+t}{t},
\ee
with $\a_0 = (1-\delta)(\max\{K+C,1\})/(2\delta)$.

Now, let us suppose that \eqref{evol_NAmt} holds for values of $m$ between $0$ and $m-1$. Let us show that it is true also for $m$. We start with the well known formulae (cf. \cite{Hu86}):\begin{multline}
\(\parcial{}{t} - \Delta\) |\nabla^{m} A |^2 \leq    - 2 |\nabla^{m+1} A|^2 \\ +  D_m  \( \sum_{i+j+k=m} |\nabla^{i} A | |\nabla^{j} A | |\nabla^{k} A | |\nabla^{m} A | + \sum_{i+j=m} |\nabla^{i} A | |\nabla^{j} \oR |  |\nabla^{m} A |  +  |\nabla^{m} \oN \oR |  |\nabla^{m} A | \), \lb{dtNAm}
\end{multline}
where $D$ is a constant which depends only on $m$ and $n$.
Also, by repetitive use of the Gauss formula for a submanifold, for the restriction of any tensor field $\overline B$ on $\oM$ to $M_t$ one has
\be \lb{GaB}
\nabla^k \overline B = \oN^k \overline B + \sum_{j=1}^k \sum_{i_0+i_1+ ... +i_j = k-j} \oN^{i_0}\overline B * \nabla^{i_1}A * \dots * \nabla^{i_j}A,
\ee
where \lq\lq $\ *$ " has the same meaning that in \cite{Hu86}. 

If, for every $k$, $|\oN^k\overline B|$ is bounded in $\oM$ by some constant $\overline b_k$, from  \eqref{GaB} and the induction hypothesis we obtain, renaming the constants each time we need, 
\begin{align}
|\nabla^k \overline B|  
& \le \overline b_k +  c(n,k) \sum_{j=1}^k  \sum_{i_1+ ... +i_j \le k-j}  \overline b_{k-j-(i_1+ ... +i_j)}\(\a_{i_1} \dots  \a_{i_j}\)^{1/2} \(\fracc{1+t}{t}\)^{(i_1+\dots+i_j+j)/2} \nn\\
& \le \overline b_k +  \overline c(n,k,\oM,T) \(\fracc{1+t}{t}\)^{k/2} \le b_k \(\fracc{1+t}{t}\)^{k/2} , \lb{b_NkB}
\end{align}
where $b_k = \max\{\overline b_k, \overline c(n,k,\oM,T) \}$. 

Like in \cite{EcHu0}, we consider the function $f_m= t^{m+1} |\nabla^m A|^2+ \xi t^{m} |\nabla^{m-1}A|^2$, where $\xi$ is some constant that will be defined later,  and estimate its evolution with time under \eqref{vpmf}. Using \eqref{dtNAm}, we get
\begin{align}
\(\parcial{}{t} - \Delta\) f_m  & \leq   (m+1) t^m |\nabla^m A|^2 +
\\  &  + t^{m+1} \(- 2 |\nabla^{m+1} A|^2 + D_m  \( \sum_{i+j+k=m} |\nabla^{i} A | |\nabla^{j} A | |\nabla^{k} A | |\nabla^{m} A | \right. \right.\\
& \left. \left.+ \sum_{i+j=m} |\nabla^{i} A | |\nabla^{j} \oR |  |\nabla^{m} A |  +  |\nabla^{m} \oN \oR |  |\nabla^{m} A | \) \) +   m\  t^{m-1} \xi |\nabla^{m-1} A|^2\\
 &  + t^{m} \xi  \(- 2 |\nabla^{m} A|^2 +  D_{m-1}  \( \sum_{i+j+k=m-1} |\nabla^{i} A | |\nabla^{j} A | |\nabla^{k} A | |\nabla^{m-1} A | \right.\right. \\
 & \left.\left. + \sum_{i+j=m-1} |\nabla^{i} A | |\nabla^{j} \oR |  |\nabla^{m-1} A |  +  |\nabla^{m-1} \oN \oR |  |\nabla^{m-1} A | \)\)
\lb{dtNfm}
\end{align}
From  the induction hypothesis, and because $0 \le (t^{k/2}  |\nabla^{k} A | - t^{m/2}  |\nabla^{m} A |)^2 = t^{k}  |\nabla^{k} A |^2  + t^{m}  |\nabla^{m} A |^2 - 2 t^{k/2}  |\nabla^{k} A |  t^{m/2}  |\nabla^{m} A | $,
\begin{align}\nn
t^{m+1} & \sum_{i+j+k=m} |\nabla^{i} A | |\nabla^{j} A | |\nabla^{k} A | |\nabla^{m} A |   \le t^{m+1} \sum_{i+j+k=m} (\a_i \a_j)^{1/2} \fracc{(1+t)^{\frac{i+j}{2}+1}}{t^{\frac{i+j}{2}+1}} |\nabla^{k} A | |\nabla^{m} A |  \end{align}
\begin{align}
& = c \sum_{k=0}^m (1+t)^{\frac{m-k}{2}+1} t^{k/2} \ |\nabla^{k} A | t^{m/2}  |\nabla^{m} A |  \nn
\end{align}
\begin{align}
& \le (c/2) (1+t)^{m+1} \sum_{k=0}^m  \(\fracc{t^{k}}{(1+t)^k}  |\nabla^{k} A |^2 + \fracc{t^{m}}{(1+t)^m}  |\nabla^{m} A |^2\) \nn \\
& = (c/2) (1+t)^{m+1}  \(\sum_{k=0}^{m} \fracc{t^{k}}{(1+t)^k}   |\nabla^{k} A |^2 +(m+1)  \fracc{t^{m}}{(1+t)^m}  |\nabla^{m} A |^2 \) \nn\\
& \le (m+2) (c/2) (1+t)^{m+1}  \sum_{k=0}^{m} \fracc{t^{k}}{(1+t)^k}  |\nabla^{k} A |^2 . \le c_m (1+t)^{m+1} \sum_{k=0}^{m}  \fracc{t^{k}}{(1+t)^k}  |\nabla^{k} A |^2 \lb{Tijkm}
\end{align}

For the second summand in \eqref{dtNfm}, denoting by  $r_j$ the analog of $b_j$ in \eqref{b_NkB} when $\overline B$ is $\oR$, using \eqref{b_NkB} and the induction hypothesis
\begin{align}
t^{m+1} \sum_{i+j=m}   |\nabla^{i} A | |\nabla^{j} \oR |  |\nabla^{m} A | & \le  t^{m+1} \sum_{i+j=m} \a_i^{1/2} \fracc{(1+t)^{(i+j+1)/2}}{t^{(i+j+1)/2}} r_j  |\nabla^{m} A | \nn \\
 &\le  (1+t)^{m+1}   \( \sum_{i+j=m} \a_i^{1/2} r_j\) \(1+\fracc{t^{m+1}}{(1+t)^{m+1}}   |\nabla^{m} A |^2\) \nn \\
 & =  (1+t)^{m+1}   d_m + d_m  t^{m+1} |\nabla^{m} A |^2. \lb{Tijm}
 \end{align}
Analogously,   denoting by  $\tilde r_j$ the analog of $b_j$ in \eqref{b_NkB} when $\overline B$ is $\oN\oR$, 
\begin{align}
t^{m+1} |\nabla^{m} \oN \oR |  |\nabla^{m} A | & \le (1+t)^{m/2} \ t^{(m/2)+1}\  \tilde r_{m} | \nabla^m A| =  (1+t)^{m}\  t\   \tilde r_m \  \fracc{t^{m/2}}{(1+t)^{m/2}} | \nabla^m A| \nn \\
& \le (1+t)^{m}\  t\   \tilde r_m \(1+\fracc{t^{m}}{(1+t)^{m}} | \nabla^m A|^2\)  \lb{Tmm}
\end{align}
and
\be\lb{Tijkm1}
t^{m}  \sum_{i+j+k=m-1} |\nabla^{i} A | |\nabla^{j} A | |\nabla^{k} A | |\nabla^{m-1} A |   \le c_{m-1} (1+t)^{m} \sum_{k=0}^{m-1}  \fracc{t^{k}}{(1+t)^k}  |\nabla^{k} A |^2 
\ee
\begin{align}
t^{m} \sum_{i+j=m-1}   |\nabla^{i} A | |\nabla^{j} \oR |  |\nabla^{m-1} A | & \le  (1+t)^{m}   d_{m-1} + d_{m-1} (1+t)^{m} \fracc{t^{m}}{(1+t)^{m}} |\nabla^{m-1} A |^2 \lb{Tijm1}
 \end{align}
 \begin{align}
t^{m} |\nabla^{m-1} \oN \oR |  |\nabla^{m-1} A | \le (1+t)^{m-1}\  t\  \tilde r_{m-1}\  \( 1 + \fracc{t^{m-1}}{(1+t)^{m-1}}| \nabla^{m-1} A|^2\)  \lb{Tmm1}
\end{align}
By substitution of all these inequalitues  in \eqref{dtNfm} we get
\begin{align}
\(\parcial{}{t} - \Delta\) f_m  & \leq   (m+1) t^m |\nabla^m A|^2 
 +   D_m  c_m (1+t)^{m+1} \sum_{k=0}^{m}  \fracc{t^{k}}{(1+t)^k}  |\nabla^{k} A |^2  \nn \\
&  + D_m  (1+t)^{m+1}   d_m+ D_m d_m    t^{m+1} |\nabla^{m} A |^2 \nn  \\  
& +  D_m (1+t)^{m} \ t\  \tilde r_{m}  \(1 + \fracc{t^m}{(1+t)^m} |\nabla^m A| ^2\) \nn \\
 &+   \xi\ m\  t^{m-1} |\nabla^{m-1}A|^2 \nn\\
 &  - 2 t^{m} \xi   |\nabla^{m} A|^2 +  D_{m-1} \xi  c_{m-1} (1+t)^{m} \sum_{k=0}^{m-1} \fracc{t^{k}}{(1+t)^k}  |\nabla^{k} A |^2  \nn \\
 & + \xi D_{m-1}  (1+t)^{m}   d_{m-1} +  \xi D_{m-1} d_{m-1} (1+t)^{m} \fracc{t^m}{(1+t)^m} |\nabla^{m-1} A |^2 \nn \\
 &+ \xi D_{m-1}(1+t)^{m-1} \ t \ \tilde r_{m-1} \(1 + \fracc{t^{m-1}}{(1+t)^{m-1}} | \nabla^{m-1} A|^2\) 
\lb{dtNfmb}
\end{align}
Using again the induction hypothesis, grouping terms and renaming constants,
\begin{align}
\(\parcial{}{t} - \Delta\) f_m  & \leq   \( (m+1) \fracc{1}{t} +    D_m  c_m \fracc{1+t}{t} + D_m d_m   + D_m  \  \tilde r_{m}  -  \ \fracc{2}{t} \xi    \)  t^{m+1}\  |\nabla^{m} A |^2  \nn \\
& +   D_m  c_m \fracc{(1+t)^{m+2}}{t}  \sum_{k=0}^{m-1}   \a_k   + D_m  (1+t)^{m+1}   d_m 
+  D_m (1+t)^{m} \ t\  \tilde r_{m}   \nn \\  
& +   \xi\ m\    \fracc{(1+t)^m}{t}  \a_{m-1}  +  D_{m-1} \xi  c_{m-1} \fracc{(1+t)^{m+1}}{t}   \sum_{k=0}^{m-1}  \a_k  \nn \\
 & + \xi D_{m-1}  (1+t)^{m} \  d_{m-1} +  \xi\  D_{m-1} d_{m-1} (1+t)^{m} \ \a_{m-1}\nn \\
 &+ \xi D_{m-1}(1+t)^{m-1} \ t \ \tilde r_{m-1}  + \xi D_{m-1}(1+t)^{m}  \ \tilde r_{m-1}\a_{m-1}  
\lb{dtNfmd}
\end{align}
Let us rename $D_m c_m = C_m$, $D_m(d_m+\tilde r_m)= E_m$, and let us define $\xi$ by
\be\lb{defxi}
2 \xi = m+1 + C_m  (1+2T) + E_m 2 T
\ee
With this choosing of $\xi$,  the coefficient of  $t^{m+1}  |\nabla^{m} A |^2$ in \eqref{dtNfmd} is bounded from above by 
\begin{align}
 \( C_m\fracc{t-2T}{t} + E_m \(1- 2 \fracc{T}{t}\)\) = F_m \fracc{t-2T}{t}  <0. \lb{ineqxi}
\end{align}

From \eqref{dtNAm} to \eqref{ineqxi} the computations are valid for compact and non-compact cases. When $M$ is compact,
let $t\in ]0,T[$. Let $t_0\in [0,t]$ where $f_m$ attains its maximum value. At the point $(p,t_0)$ where this happens, it follows from \eqref{dtNfmd}, \eqref{defxi},  \eqref{ineqxi} and the maximum principle, renaming again all the constants, that
\begin{align}
t_0^{m+1}  |\nabla^{m} A |^2 & \leq  
\( \mathcal A  \fracc{(1+t_0)}{2T-t_0}  + \mathcal B  \fracc{t_0}{2T-t_0}+ \mathcal C \fracc{t_0^2}{(1+t_0)(2T-t_0)}  \right. \nn \\
 & \qquad + \mathcal D (2 T-t_0) + \mathcal E \fracc{t_0}{(1+t_0)(2T-t_0)} +\mathcal  F  \fracc{1}{(1+t_0)(2T-t_0)}\nn \\ & \qquad \left. +  \mathcal G \fracc{t_0^2}{(1+t_0)^2(2T-t_0)}  \) (1+t_0)^{m+1}
\lb{dtNfmf}
\end{align}

Then, by the definition of maximum and because $t_0\le t < T$
\begin{align}
t^{m+1} |\nabla^{m} A |^2 & \le f_m (t_0) = t_0^{m+1} |\nabla^m A|^2+ \xi t_0^{m} |\nabla^{m-1}A|^2 \nn \\
& \le \( C_1 T + C_2 + C_3 \frac1T\) (1+t)^{m+1}
 + K (1+t)^{m+1} \nn \\
 & \le \mathcal K (1+t)^{m+1}.
\lb{dtNfmg}
\end{align}
which finishes the proof by induction of the  formula \eqref{evol_NAmt} when $M$ is compact.

If $M$ is non-compact,  we start using the inequality \eqref{phg2} to obtain the following inequation for the evolution of $t \phi \frak g$
\begin{align}\lb{evol_tpg}
\(\parcial{}{t} - \Delta\) (t\phi  \frak g) & \le  - \<\nabla(t \phi\g),   \fracc{\nabla\psi}{\psi} + 2 \frac{\nabla\phi}{\phi}\> - 2 \delta \ t\ \phi\ \g^2 + \fracc{\nabla \phi}{\sqrt{\psi}}\ t\  \g \sqrt{\g}  \nn \\
&\qquad +  \( 6 \a  \cm(R)^2 + \eta |\nabla\phi|  + K \phi\) t\ \g + C\ t\  \sqrt{\phi}\sqrt{\phi \g} + \phi \ \frak g.
\end{align}
Given $t \in ]0,T[$, let $t_0\in[0,t]$ be the time when $t \phi \frak g$ attains it maximum value $t_0 \phi  \frak g_0$ in $S_{\r_1,t}$ (then $\phi \frak g_0 = \max_{x\in S_{\r_1,t_0}} \phi \frak g(x)$ and $t_0\ne 0$). By the maximum principle, from \eqref{evol_tpg} and, multiplying by $\fracc{\sqrt{t_0\phi}}{\sqrt{\g_0} }$, we get 
\begin{align}
2 \delta \(\sqrt{ t_0 \phi \g_0 }\)^3 \le \fracc{\nabla \phi}{\sqrt{\psi}\sqrt{\phi}} \sqrt{t_0}\  t_0 \phi\g_0 + &  \(6 \a  \cm(R)^2 + \eta |\nabla\phi|   + K \phi \) {t_0} \sqrt{t_0 \phi \g_0}\nn \\
  &+ C \(\sqrt{t_0}\)^3 \(\sqrt{\phi}\)^3 + \phi \sqrt{t_0 \phi g_0}. \lb{phg5}
\end{align}
 Having into account \eqref{insqph}, that $\sm(\r_1)\le \cm(\r_1)$ and dividing by $\cm(\r_1)^3$, we get
\begin{multline}
2 \delta \( \fracc{\sqrt{t_0 \phi\ \g_0}}{\cm(\r_1)}\)^3 - 2 \sqrt{\a\ T}  \sqrt{1-\delta}   \( \fracc{\sqrt{t_0\phi\ \g_0}}{\cm(\r_1)}\)^2  
\\ - \( \( \( 6 +\fracc{K}{\mu^2} \) \a\  +  \fracc{\eta 2 \a}{-\mu} \)  T + \fracc{\a^{3/2} }{-\mu^3} \) \fracc{\sqrt{t_0\phi\ \g_0}}{\cm(\r_1)}   - C \a^{3/2} \fracc{1}{-\mu^3} T \le 0.  \lb{phg4t} 
\end{multline}
From now, arguing like after \eqref{phg4}, since now $t\phi \g$ is $0$ at $S_{\r_1,0}$,
\be 
 \( \fracc{ \sqrt{t \phi\ \g}}{\cm(\r_1)}\)  \le  \( \fracc{ \sqrt{t_0\phi\ \g_0}}{\cm(\r_1)}\) \le D, \quad \text{ and } \quad  t  |A|^2 \le \fracc{ \mu^2 (1-\delta)  D^2}{\a (1-\gamma)^2}, \label{oot}
\ee 
   This shows that, on $ \mathbb S_{\r_1,T,\gamma}$, $t |A|$ is bounded by a bound depending   $T$ and  $\r_1$ (only through $\delta$) and on  $M_0$ (only through the maximum of $v$ in $\mathbb S_{\r_1,T}$). Renaming the constants, and, since $1 \le 1+t$,
   $$  t  |A|^2 \le \a_0 (1+t) \text{ on } \mathbb S_{\r_1,T,\gamma}= \mathbb S_{\r_2,T}.$$
     This proves \eqref{evol_NAmt} for $M$ non-compact and $m=0$. To prove it for every $m$ we consider the evolution of $\phi f_m$. From \eqref{evol_phi} and \eqref{dtNfmd}, computing like we did for getting \eqref{phg0} and renaming some constants,
     \begin{align}
\(\parcial{}{t} - \Delta\) (\phi f_m)  & \le \phi \(   \(  \fracc{(m+1)}{t} +    C_m \fracc{1+t}{t} + E_m   -  \fracc{2}{t} \xi    \) {t^{m+1}}  |\nabla^{m} A |^2  \right. \nn \\
& +   \mathcal C \fracc{(1+t)^{m+2}}{t}    + \mathcal D  (1+t)^{m+1}   + \mathcal B (1+t)^{m} \ t\    \nn  \\  
& \left. +   \xi\  \(\mathcal E\    \fracc{(1+t)^m}{t}   +  \mathcal F \fracc{(1+t)^{m+1}}{t}     + \mathcal G  (1+t)^{m}  
 + \mathcal H (1+t)^{m-1} \ t \)\) \nn \\
 & -2  \a\ f_m \frac{|\nabla \phi|^2}{ \phi'^2} - 2 \<\frac{\nabla \phi}{\phi},\nabla ( \phi  f_m) - f_m \nabla\phi\>  .
 \end{align}

From here, reasoning like we did from equations \eqref{phg0f} to \eqref{oom}, but taking, at the end,  $\xi = m+1+C_m(1+2 T) + E_m 2 T + 3  \fracc{\mu^2\ T}{(1-\gamma)^2}$ (instead of  \eqref{ineqxi}), we get $t^{m+1}|\nabla^m A|^2 \le \a_m (1+t)^{m+1}$ on $\mathbb S_{\r_{1+m},T,\gamma}$, which finishes the proof by induction for the non-compact case.
\end{demo}

\begin{teor} \label{steL} 
Let $M$ be complete (compact or not).  If $\oM = M \times_{\p}\re$ and $M_0$ is  a Lipschitz  continuous graph over $M$, then there is a  solution $M_t$ of  \eqref{vpmf}   with initial condition $M_0$ which is a graph over $M$, is defined on $[0,\infty[$, and $M_t$ is smooth for every $t\in]0,\infty[$. 
\end{teor}
\begin{demo}
We shall write the details for the case $M$ non-compact. When $M$ is compact the arguments are similar, simplified by the fact that we can take always $M$ or $M_t$ instead of $\tB(x_0,\r_0)$ or $\mathbb S_{\r_k,t,\gamma}$ respectively. 

Let $M_{k}$ be a sequence of smooth manifolds given by smooth graphs $u_k$ over $M$ and converging to $M_0$. For each $M_k$, let $M_{kt}$ be the smooth solution  of \eqref{mcf} which has $M_k$ as initial condition, which exists and is defined for $t\in[0,\infty[$ by theorems \ref{existco} and \ref{stenonco}. Each $M_{kt}$ is represented by the graph of a function $u_k(\cdot, t)$ which is a solution of \eqref{mcft} with the initial condition $u_k(\cdot, 0)$. 

Given $\r_0>0$ and $T>0$, let us define $\r$ by \eqref{R'R0}.  It follows from \eqref{vboundex} and \eqref{nuv} that 
\be
\lb{boundtNu} |\tN u_{k}| \le c \text{ on } \mathbb S^k_{\rho,T,\gamma} = \{ (x,u_k(x,t)); \fracc{\cm(\ttr((x,u_k(x,t)))}{- \mu} \le \gamma \fracc{\cm(\r)}{- \mu} e^{-\beta n t}, t\in [0,T]\}, 
\ee
 where $c$ is a constant which depends only on $\r_0$ (through $\r$), $\gamma$, $T$,  an upper bound   of $v$ on $\mathbb S^k_{\r,0}\subset M_k$ and the bounds of $\p$. Since $M_0$ is Lipschitz, $v$ is bounded on the corresponding $\mathbb S_{\r,0}\subset M_0$. Since the $M_k$ converge to $M_0$, the maxima  of $v$ on $\mathbb S^k_{\r,0}\subset M_k$ will be near the maximum of $v$ on $\mathbb S_{\r,0}\subset M_0$ for $k$ big enough. Then we can take the constant $c$ in \eqref{boundtNu} independent of $k$. 
From the definition of $\r$, if $(x,t)\in \tB(x_0,\r_0) \times [0,T]$ then $(x,u_k(x,t)) \in \mathbb S_{\r,T,\gamma} $ and \eqref{boundtNu} holds on $\tB(x_0,R) \times [0,T]$.
 
 The same argument as given at the end of last paragraph shows that, by Lemma \ref{boundNAmL}, 
\be\lb{boundNAmLR}
 |\nabla^m A|^2 \le \alpha_m \(\fracc{t+1}{t}\)^{m+1} \text{ on }  \tB(x_0,\r_0) \times ]0,T],
 \ee
where the constants $\a_m$ depend only on $m$, $\r_0$, $\gamma$ and $T$. 

Now, let us consider the functions $u_{k\r_0}(x,t)= u_k|_{\tB(x_0,\r_0) \times [0,T]}(x,t)$.
Through the relation between $u$ and $A$, the bounds \eqref{boundtNu} and \eqref{boundNAmLR} give that for every $t\in]0,T]$,  there are constants $\beta(m,t,\r_0,T,\gamma)$  depending only on $m$, $t$, $\r_0$, $T$ and $\gamma$ such that
\be
|\tN^m u_{k\r_0}(\cdot, t)| \le \beta(m,t,\r_0,T,\gamma).
\ee 
then, by the Ascoli-Arzela lemma, for each  $t$ and each $\r_0$, there is a subsequence $u_{k\r_0}(\cdot, t)$ which converges to a $C^\infty$ function $u_{0\r_0}(\cdot, t)$. 

If we consider the $u_{k\r_0}(\cdot, t)$ and $u_{0\r_0}(\cdot, t)$ as functions of $t$, the same arguments used in section \ref{STE} (now with bound \eqref{evol_NAmt}) show that  $\ds \left|\parcial{u}{t}\right|$ and $\ds\left|\parcial{^2u}{t^2}\right|$ are bounded on $\tB(x_0,\r_0)) \times[t_0,T]$ by a constant depending only on $t_0>0$, $\r_0$, $T$ and $\gamma$, but not on $k$. 
Once we know that, we can conclude that the limit $u=u_{0\r_0}$ is, at least,  of class $C^1$ on $t$ and, since all $u_{k\r_0}$ satisfy equation \eqref{mcfg} on $\tB(x_0,\r_0)) \times[t_0,T]$, also does $u_{0\r_0}$. Taking $\r_0$ and $T$ bigger, we have sequences of functions satisfying similar conditions and which coincide with the older ones for the older values of $\r_0$ and $T$, then for these new values of $\r_0$ and $T$ we have a new limit function which coincides withe the older limit when restricted to the older values of $\r_0$ and $T$. Letting $\r_0$ an $T$ go to infinity and $t_0$ to $0$, this gives  a function $u_0$ which is a solution of \eqref{mcfg}, satisfies the initial condition $u_0(\cdot,0)=$ the function defining $M_0$ as a graph, and, for every $t\in]0,\infty[$, $u_0(\cdot,t)$ is $C^\infty$ as a function of $M$.
\end{demo}


\section{Some results on convergence}\lb{CONV}

\begin{lema}\label{ltozero}
Let $\oM$ be a complete (may be compact) riemannian manifold, with bounded sectional curvature  $k_0 \le \oSec \le k <0$ and having a complete totally geodesic hypersurface $M$. Let $M_t$ be the evolution of a complete hypersurface $M_0$ at time $t$ by  \eqref{vpmf}. Let us suppose that, for every $x\in M_0$, the distance $\ell(x)$ from $x$ to $M$ is bounded from above by some constant $\ell_0$. If $M$ is non-compact, we add the hypothesis that the norm $|A_t|$ of the Weingarten map $A_t$  of the hypersurface $M_t$  is bounded by a constant (depending on $t$). Then,  one has
\be \label{bol}
\sk(\ell(F(x,t)))  \le \sk(\ell_0) e^{k n t}
\ee
\end{lema}
\begin{demo}
A calculation similar to that done before for $u$, but now for $\ell$ is
\begin{align}
(\nabla^2 \ell)(E_i,E_j) & = \< \nabla_{E_i} \nabla \ell, E_j\> =  \< \nabla_{E_i} (\partial_\ell -\<\partial_\ell, N\> N), E_j\>\nn \\
&=  \< \oN_{E_i}   \partial_\ell, E_j\> + \<\partial_\ell, N\> h(E_i,E_j) \nn\\
&=     \<\oN_{E_i - \<E_i, \partial_\ell\>\partial_\ell} \partial_\ell, E_j\> + \<\partial_\ell, N\> h(E_i,E_j) 
\end{align}
and, using \eqref{comptg}, one gets
\begin{align}
\Delta \ell &=   \sum_i\<\oN_{E_i - \<E_i, \partial_\ell\>\partial_\ell} \partial_\ell, {E_i - \<E_i, \partial_ell\>}\partial_\ell\> + \<\partial_\ell, N\>H\nn\\
& \ge - k \fracc{\sk}{\ck} \(n - |\partial_\ell^\top|^2\) +  \<\partial_\ell, N\> H  .\label{Dell}
\end{align}
Then, for the evolution of $\sk(\ell):=\sk\circ\ell$ one has
\begin{align}
\parcial{\sk(\ell)}{t} &= \ck(\ell)\ H \<\partial_\ell,N\> \le \ck(\ell)\  \Delta \ell + k \sk(\ell) \(n - |\partial_\ell^\top|^2\) \nn \\
&= \Delta \sk(\ell) - (-k \sk(\ell) |\partial_\ell^\top|^2)+ k \ \sk(\ell) \(n - |\partial_\ell^\top|^2\) \nn \\
& =\Delta \sk(\ell) + n\ k \ \sk(\ell) \lb{evolskl}
\end{align}
If $M$ is compact, then the maximum principle states that $\ell$ is bounded by the solution of the ODE $\sk(\ell)' = n\ k \ \sk(\ell) $ with the initial condition $\ell(0)=\ell_0$, which is $\sk(\ell) = \sk(\ell_0) e^{n k t}$, from which the statement of the theorem follows.

If $M$ is non-compact, let us consider first the evolution of $\ell$, which follows from \eqref{Dell} as above
\begin{align}
\parcial{\ell}{t} &=  H \<\partial_\ell,N\> \le  \Delta \ell +  k \fracc{\sk(\ell)}{\ck(\ell)}  \(n - |\partial_\ell^\top|^2\) \le \Delta \ell +  k \fracc{\sk(\ell)}{\ck(\ell)}  (n-1) \le \Delta \ell, \lb{evoll}
\end{align}
we observe first that the added hypothesis implies that we can apply the maximum principle for noncompact manifolds given in Lemma \ref{echu43} to $f=\ell-\ell_0$. 
In fact,  our hypothesis, together with the hypothesis $k_0 \le \oSec$ and the Gauss formulae give that the sectional curvature of each $M_t$ is bounded from below, which shows that the hypothesis on the growing of volume of geodesic balls of that Lemma is satisfied. Conditions (i) and (ii) are obvious. Condition iv) follows from the evolution equation $\parcial{}{t}g = - 2 H \a$ and the added hypothesis. Condition (iii) follows from $|\nabla \ell | \le 1$ and the fact that  $d\mu_t \le \fracc{\s^n(r)}{r^n} d\mu_e$, where $\lambda$ is determined by $\mu_0$ and the upper bounds of $|A_t|$ and $d\mu_e$ is the euclidean measure. Then we can conclude, from the quoted theorem, that $\ell- \ell_0 \le 0$ for all time.

Once we know that $\ell$ is bounded, we can study the evolution of $u=\sk(\ell) - \sk(\ell_0) e^{nkt}$, which, from \eqref{evolskl} is
\begin{align}
\parcial{u}{t} &=\Delta u + n\ k \ u, \lb{evolu}
\end{align}
and, again, this equation satisfies the conditions of Theorem 4.3 in \cite{EcHu}, the only condition that has to be checked now is (iii), which follows because $|\nabla u| = \ck(\ell) |\nabla\ell| \le \ck(\ell_0)$ and now we know that $\ell \le \ell_0$. Then, as before, we can conclude that $u\le 0$ and the theorem is proved.
\end{demo}

With the same notations \eqref{def_mueta} and \eqref{def_nu}, we have
\begin{teor}\label{convcomp} 
Let $M$ be compact. Let $M_t$ be the solution of  \eqref{vpmf}  defined on a maximal time interval $[0,T[$. If $\oM = M \times_{\p}\re$,      $0 > -\mu_1 > \tSec \ge \fracc{- n \mu_1 + \mu_2}{n-1} $ and   $M_0$ is a graph over $M$,   then $M_t$ is a graph, $T=\infty$ and $M_t$ converges in the  $C^\infty$ topology to  $M\times \{0\}$ as $t\to\infty$. 
\end{teor}
\begin{demo}  
Notice that the hypothesis on the lower bound of $\tSec$ is equivalent to $\nu\le 0$, then, from theorems \ref{graphbo2comp} and \ref{existco} it follows that both $v$ and $|\nabla^iA|$, $i=0,1,...$ are bounded by some constants not depending on time, then, by Ascoli-Arzelˆ, the $M_t$ converge to some limit $M_\infty$. The hypotheses  $-\mu_1 > \tSec$ implies $-\mu_1 > \oSec$, then we can apply  Lemma \ref{ltozero} to conclude that $M_\infty$ must be at distance $0$ from $M$, then it must be $M$.
\end{demo}
\begin{teor}\label{convnonco} 
Let $M$ be complete non compact with a pole $x_0$.  Let $M_t$ be the solution of  \eqref{vpmf}  defined on a maximal time interval $[0,T[$. If  $\oM = M \times_{\p}\re$,     $0 > -\mu_1 > \tSec \ge \fracc{- n \mu_1 + \mu_2}{n-1} $,   $M_0$ is a graph over $M$ with $|A|^2$ bounded,   then $M_t$ is a graph, $T=\infty$ and $M_t$ converges in the  $C^\infty$ topology to  $M\times \{0\}$ as $t\to\infty$. 
\end{teor}
\begin{demo}  
Looking  at \eqref{vboundnoR} and at the proof of Theorem   \ref{existnonco}, one observes that, when  $v$ and $|\nabla^iA|$, $i=0,1,...$, are bounded on $M_0$, then both $v$ and $|\nabla^iA|$ are bounded by some constants depending on time, but not on $R$. Then, we can apply Lemma \ref{echu43} to the evolution equation \eqref{dvtnb} satisfied by $v$ when $\nu\le 0$ to conclude that $v$ is bounded by a bound not depending on time. Using this fact in the proof of Theorem \ref{existnonco}, we see that the bounds of $|\nabla^iA|$ do not depend on time, then, as above, we can apply Ascoli-Arzela to conclude that $M_t$ has a limit $M_\infty$, and, using   Lemma \ref{ltozero} as before, we conclude that $M_\infty$ must be  $M$.
\end{demo}

\section{Appendix: About the flow in $\re\times_\p M$}

When we consider the ambient space $\oM= \re\times_\p M$, that is $(\re\times M, \og= du^2 + \p(u) \tg)$, where $\p:\re \flecha \re$. If we consider the evolution under \eqref{vpmf} of a hypersurface $M_0$ which is a graph over $M$ we obtain for $v=\<N,\partial_y\>^{-1}$ the evolution equation
 \begin{align}
\parcial{v}{t} & = \Delta v - \frac2{v} |\nabla v|^2 +  2 \fpp  \<\nabla v , \nabla u\> - \fpp 2 H v^2    \nn\\
& \qquad -  \(1- \fracc1{v^2}\)   \(\fracc{\tRic_{\t1\t1}}{\p^2}  + n  \fppp   \) v
-   |A|^2  v   -  \(\fpp\)^2 n \frac1{v}.
 \label{dvtn}
\end {align}
where $\tRic_{\t1\t1}$ is the Ricci curvature in the direction $\tN u$.

The term $- \fpp 2 H v^2$  in this equation do not allow to apply the maximum principle as we did with
\eqref{dvtnb} and, in fact, can be considered as the analytical reason why the property of being a graph is not preserved in general under the MCF in $ \re\times_\p M$. 

Looking again at \eqref{dvtn} it is possible to think that we can get interesting results for hypersurfaces satisfying that the sign of $H$ and the sign of $\p'$ is the same. But this condition has many drawbacks :

\vspace{-0.3cm}
\begin{enumerate}
\item A computation shows that if $M_0$ is the graph of a function $u$,  
$$\Delta u  = \(n-\(1-\fracc1{v^2}\)\) \fpp + \fracc1{v} H . $$
 Then, when $M$ is compact, it follows from this formula and the maximum principle that $H$ has the same sign that $\p'$ if and only if $u=0$.
\item When $M$ is non-compact, the known proof (for the euclidean space) given in \cite{EcHu} that the sign of $H$ is preserved uses the bounds of $|A|$, and, for them, bounds of $v$ are also used. But in our case, are just the bounds of $v$ which we do not know.
\end{enumerate}

\vspace{-0.2cm}
Another idea of why it is easy that the property of being a graph is preserved for $M\times_\p\re$ and not for $\re\times_\p M$ is given by the following pictures. In both $M=\mathcal H^n$ and $\p(u)=\cosh u$ for the first and $\p(x)=\cosh(\dist(x_0,x))$ for the second. Then $\oM = \mathcal H^{n+1}$ in both, but a graph for the first one corresponds to the sense of graph in this appendix (a graph for geodesics) and the second one corresponds to a graph in the sense o the previous sections of this paper, that is a graph for \lq\lq equidistant" curves. 

\begin{center}
\includegraphics[width=6cm,height=6cm]{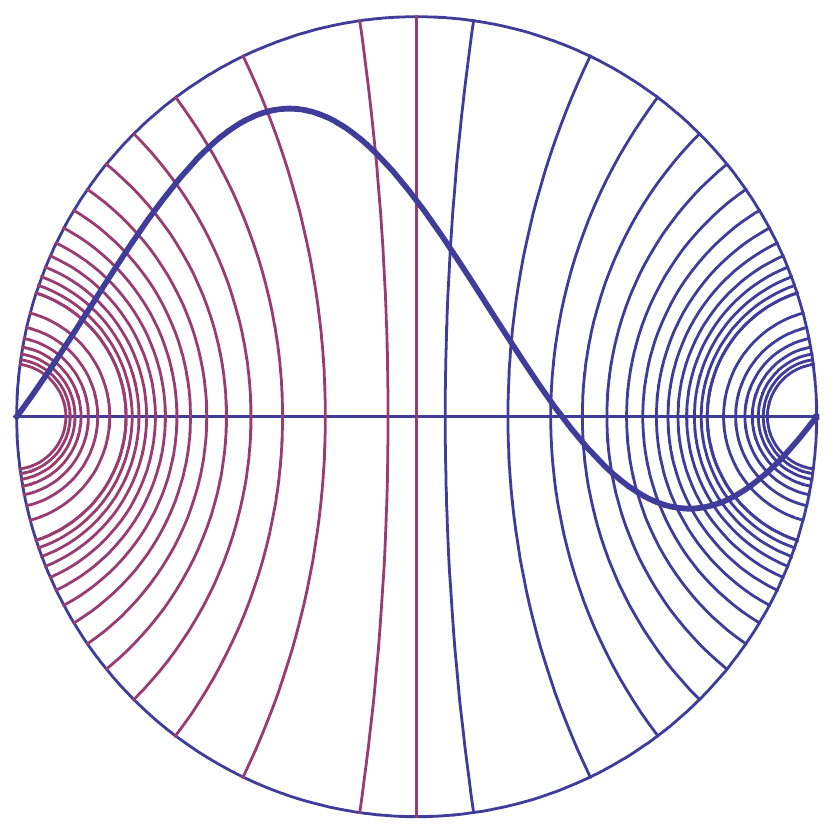}
\put(-130,-10){a graphic for geodesics }\hskip1truecm
\includegraphics[width=6cm,height=6cm]{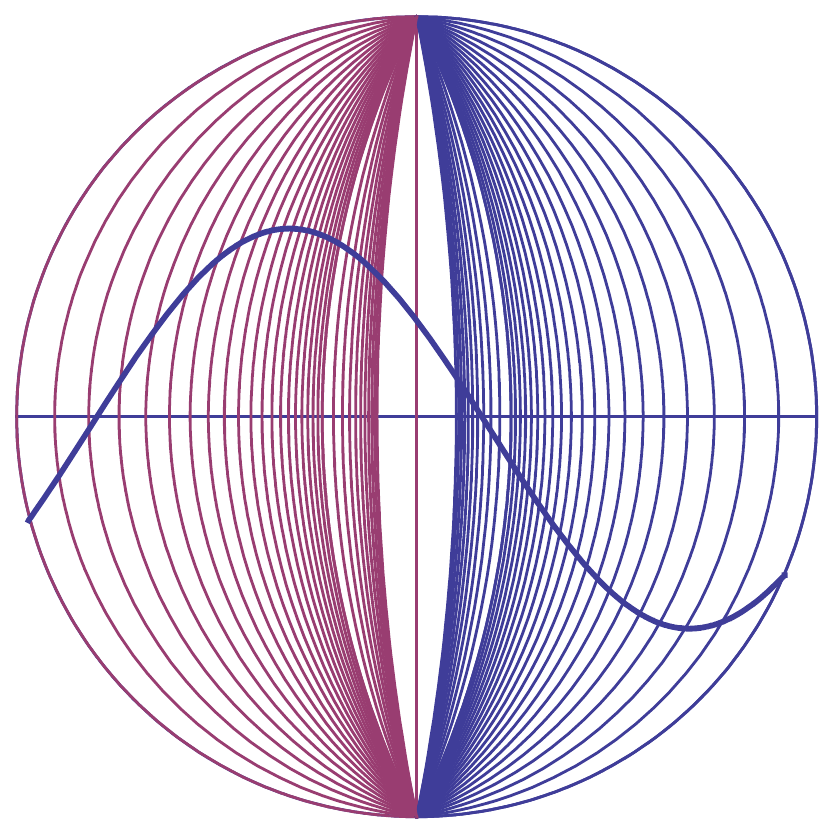}
\put(-130,-10){a graphic for equidistants  which } \put(-132,-22){is not a graphic for geodesics}

 \end{center}

{\footnotesize

\bibliographystyle{alpha}

}

\vskip1truecm

{\small 

Address

\begin{tabular}{ c c }
 Kharkov National University &Universidad de Valencia\\
Mathematics Faculty. Geometry Department &Departamento de Geometr\'{\i}a y Topolog\'{\i}a\\
 Pl. Svobodi 4 &Avda. Andr\'es Estell\'es, 1\\
 61077-Kharkov, Ukraine & 46100-Burjassot (Valencia) Spain\\
       email: borisenk@univer.kharkov.ua & email: miquel@uv.es  \\

\end{tabular}

}

\end{document}